\documentclass[12pt]{amsart}
\usepackage{amssymb}

\begin{document}
\numberwithin{equation}{section}
\setlength{\emergencystretch}{3em}

\def\Label#1{\label{#1}}

\def\1#1{\ov{#1}}
\def\2#1{\widetilde{#1}}
\def\3#1{\mathcal{#1}}
\def\4#1{\widehat{#1}}

\def\s{s}
\def\k{\kappa}
\def\ov{\overline}
\def\span{\text{\rm span}}
\def\tr{\text{\rm tr}}
\def\GL{{\sf GL}}
\def\xo {{x_0}}
\def\Rk{\text{\rm Rk\,}}
\def\sg{\sigma}
\def \emxy{E_{(M,M')}(X,Y)}
\def \semxy{\scrE_{(M,M')}(X,Y)}
\def \jkxy {J^k(X,Y)}
\def \gkxy {G^k(X,Y)}
\def \exy {E(X,Y)}
\def \sexy{\scrE(X,Y)}
\def \hn {holomorphically nondegenerate}
\def\hyp{hypersurface}
\def\prt#1{{\partial \over\partial #1}}
\def\det{{\text{\rm det\,}}}
\def\wob{{w\over B(z)}}
\def\co{\chi_1}
\def\po{p_0}
\def\fb {\bar f}
\def\gb {\bar g}
\def\Fb {\ov F}
\def\Gb {\ov G}
\def\Hb {\ov H}
\def\zb {\bar z}
\def\wb {\bar w}
\def \qb {\bar Q}
\def \t {\tau}
\def\z{\chi}
\def\w{\tau}
\def\Z{\zeta}

\def \T {\theta}
\def \Th {\Theta}
\def \L {\Lambda}
\def\b{\beta}
\def\a{\alpha}
\def\o{\omega}
\def\l{\lambda}
\def \CFT{\text{\rm CFT}}

\def \im{\text{\rm Im }}
\def \re{\text{\rm Re }}
\def \Char{\text{\rm Char }}
\def \supp{\text{\rm supp }}
\def \codim{\text{\rm codim }}
\def \Ht{\text{\rm ht }}
\def \Dt{\text{\rm dt }}
\def \hO{\widehat{\mathcal O}}
\def \cl{\text{\rm cl }}
\def \bR{\mathbb R}
\def \bC{\mathbb C}
\def \bP{\mathbb P}
\def \C{\mathbb C}
\def \bL{\mathbb L}
\def \bZ{\mathbb Z}
\def \bN{\mathbb N}
\def \scrF{\mathcal F}
\def \scrK{\mathcal K}
\def \scrM{\mathcal M}
\def \cR{\mathcal R}
\def \scrJ{\mathcal J}
\def \scrA{\mathcal A}
\def \scrO{\mathcal O}
\def \scrV{\mathcal V}
\def \scrL{\mathcal L}
\def \scrE{\mathcal E}
\def \hol{\text{\rm hol}}
\def \aut{\text{\rm aut}}
\def \Aut{\text{\rm Aut}}
\def \J{\text{\rm Jac}}
\def\jet#1#2{J^{#1}_{#2}}
\def\gp#1{G^{#1}}
\def\gpo{\gp {2k_0}_0}
\def\emmp {\scrF(M,p;M',p')}
\def\rk{\text{\rm rk}}
\def\Orb{\text{\rm Orb\,}}
\def\Exp{\text{\rm Exp\,}}
\def\ess{\text{\rm Ess\,}}
\def\mult{\text{\rm mult\,}}
\def\Jac{\text{\rm Jac\,}}
\def\Span{\text{\rm span\,}}
\def\d{\partial}
\def\D{\3J}
\def\pr{{\rm pr}}
\def\dbl{[\![}
\def\dbr{]\!]}
\def\nl{|\!|}
\def\nr{|\!|}

\def \D{\text{\rm Der}\,}
\def \Rk{\text{\rm Rk}\,}
\def \ima{\text{\rm im}\,}
\def \vfi{\varphi}

\def\sideremark#1{\ifvmode\leavevmode\fi\vadjust{
\vbox to0pt{\hbox to 0pt{\hskip\hsize\hskip1em
\vbox{\hsize3cm\tiny\raggedright\pretolerance10000
\noindent #1\hfill}\hss}\vbox to8pt{\vfil}\vss}}}

\title[A Cauchy-Kowalevsky theorem for overdetermined systems]
{A Cauchy-Kowalevsky theorem for overdetermined systems of nonlinear
partial differential equations and geometric applications}
\author[M. S. Baouendi, P. Ebenfelt, D. Zaitsev]
{M. S. Baouendi, P. Ebenfelt, D. Zaitsev} \footnotetext{{\rm The
first two authors are supported in part by DMS-0701070 and DMS-0701121. The
third author is supported in part by the Science Foundation Ireland.\newline}}
\address{M. S. Baouendi, P. Ebenfelt: Department of Mathematics, University of
California at San Diego, La Jolla, CA 92093-0112, USA}
\email{sbaouendi@ucsd.edu, pebenfel@math.ucsd.edu}
\address{D. Zaitsev: School of Mathematics, Trinity College, Dublin 2, Ireland}
\email{zaitsev@maths.tcd.ie}



\newtheorem{Thm}{Theorem}[section]
\newtheorem{Def}[Thm]{Definition}
\newtheorem{Cor}[Thm]{Corollary}
\newtheorem{Pro}[Thm]{Proposition}
\newtheorem{Lem}[Thm]{Lemma}
\newtheorem{Rem}[Thm]{Remark}
\newtheorem{Ex}[Thm]{Example}
\newtheorem{Con}[Thm]{Conjecture}

\maketitle
\tableofcontents

\section{Introduction}

The main motivation for the work presented in this paper is to
construct real hypersurfaces in $\bC^{n+1}$ with maximal Levi number
(see below for the definition), a problem that has been open since
Levi numbers were introduced in \cite{BHR}. The examples constructed
here are tube hypersurfaces (i.e.\ of the form $M=\Sigma+i\bR^{n+1}$
where $\Sigma$ is a hypersurface in $\bR^{n+1}$). Moreover, we give a local
description of all such hypersurfaces. In order for the real
hypersurface $M$ to have the desired properties, $\Sigma$ must be a
non-cylindrical hypersurface whose Gauss map (or equivalently second
fundamental form) has rank one. To construct locally defined
hypersurfaces in $\bR^{n+1}$ with these properties, indeed to
parametrize all such, we prove an existence and uniqueness theorem
concerning a Cauchy problem for a class of overdetermined systems of
nonlinear partial differential equations in $\bR^n$. The Cauchy data
is posed on a noncharacteristic $k$-dimensional plane in $\bR^{n}$,
where the dimension $k$ depends, roughly, on how overdetermined the
system is. In the special case where the system is not
overdetermined at all, we have $k=n-1$ and our theorem reduces to
the classical Cauchy-Kowalevsky theorem.

Consider the Cauchy problem for an unknown vector-valued function $u:=(u^1,\ldots, u^m)$ given by
\begin{equation}\Label{CP0}
\left\{
\begin{aligned}
& \frac{\partial u^A}{\partial x^\alpha}=F^A_\alpha\left (x,u,\frac{\partial u}{\partial x^1},\ldots,\frac{\partial u}{\partial x^k}\right)\\
& u^A(x^1\ldots, x^k,x^{k+1}_0,\ldots, x^n_0)=a^A(x^1, \ldots, x^k)
\end{aligned}
\right.,
\quad A=1,\ldots, m,\ \alpha=k+1,\ldots, n,
\end{equation}
where $F^A_\alpha(x,p,p')$ are real-valued, real-analytic functions in an open subset $U$ of $\bR^n\times\bR^m\times \bR^{mk}$, $(x_0,p_0,p'_0)$ a fixed point in  $U$, and
the Cauchy data $a=(a^1,\ldots a^m)$ are real-analytic in a
neighborhood of $(x^1_0,\ldots, x^k_0)$ in $\bR^k$ with
\begin{equation}\Label{initial0}
\left(a(x^1_0,\ldots, x^k_0),\left(\frac{\partial a}{\partial x^1}(x^1_0,\ldots, x^k_0), \ldots,\frac{\partial a}{\partial x^k}(x^1_0,\ldots, x^k_0)\right)
\right )=(p_0,p'_0).
\end{equation}
We shall say that the Cauchy problem \eqref{CP0} is
{\it solvable at $(x_0,p_0,p'_0)$ for every initial data}
if there is a unique real-analytic
solution $u(x)$ to \eqref{CP0} in an open neighborhood of $x_0$ for every Cauchy data
$a(x^1,\ldots x^k)$ satisfying \eqref{initial0}.

In order to formulate our main result, we need to introduce one more notion. We shall say that a real-analytic $\bR^m$-valued function $u(x)$
defined in an open neighborhood of $x_0$
is an {\it approximate
solution at $x_0$} to the Cauchy problem \eqref{CP0} if
\begin{equation}\Label{CPapprox}
\left\{
\begin{aligned}
& \frac{\partial u^A}{\partial x^\alpha}=F^A_\alpha\left (x,u,\frac{\partial u}{\partial x^1},\ldots,\frac{\partial u}{\partial x^k}\right) +O\left (\sum_{i=1}^n|x^i-x^i_0|^2\right ),\\
& u^A(x^1\ldots, x^k,x^{k+1}_0,\ldots, x^n_0)=a^A(x^1, \ldots, x^k) +O\left (\sum_{j=1}^k|x^j-x^j_0|^3\right ),
\end{aligned}
\right.
\end{equation}
where $A$ and $\alpha$ are as in \eqref{CP0}.
Note that the property of $u(x)$ being an approximate solution at $x_0$
to the system \eqref{CP0} only depends on its $2$-jet at $x_0$
as well as the $2$-jet of the initial data $a(x^1,\ldots, x^k)$ at $(x^1_0,\ldots, x^k_0)$. The Cauchy problem \eqref{CP0} is said to be {\it approximately
solvable for every initial data} at $(x_0,p_0,p'_0)\in U$ if for every
initial data $a$
defined in a neighborhood of $(x^1_0,\ldots,x^k_0)$
and satisfying \eqref{initial0},
there is an
approximate solution $u$ to  \eqref{CP0} at $x_0$. Our main result is the following.

\begin{Thm}\Label{main1} Let
$F^A_\alpha=F^A_\alpha(x,p,p')$, for $A=1,\ldots, m$ and $\alpha=k+1,\ldots, n$, be real-valued, real-analytic functions defined in an open
subset $U$ of
$\bR^n\times\bR^m\times \bR^{mk}$. The following are equivalent.
\medskip

\noindent {\rm (i)} The Cauchy problem \eqref{CP0} is solvable
for every initial data at every
$(x_0,p_0,p'_0)\in U$.
\medskip

\noindent {\rm (ii)} The Cauchy problem \eqref{CP0} is approximately
solvable for every initial data at every $(x_0,p_0,p'_0)\in U$.

\medskip
\noindent {\rm (iii)} The functions $F^A_\alpha$ satisfy the
compatibility conditions
\begin{equation}\Label{comp10}
\Phi^A_{\alpha\beta}(x,p,p')=\Phi^A_{\beta\alpha}(x,p,p'),
\end{equation}
\begin{equation}\Label{comp20}
\Psi^{A\Gamma\Lambda}_{\alpha\beta C}(x,p,p')
= \Psi^{A\Gamma\Lambda}_{\beta\alpha C}(x,p,p'),
\end{equation}
for all $A, C=1,\ldots, m$, $\Lambda,\Gamma=1,\ldots, k$,  $\alpha,\beta=k+1,\ldots, n$, and $(x,p,p')\in U$,
where
\begin{equation}\Label{compfunc}
\begin{aligned}
\Phi^A_{\alpha\beta}:= & \frac{\partial F^A_{\alpha}}{\partial x^\beta}+\sum_{B=1}^m\frac{\d F^A_{\alpha}}{\d p^B}F^B_{\beta}+\sum_{\Gamma=1}^k\sum_{B=1}^m \frac{\d F^{A}_{\alpha}}{\d p^B_\Gamma}
\left (\frac{\d F^B_{\beta}}{\d x^ \Gamma}+\sum_{C=1}^m\frac{\d
F^B_{\beta}}{\d p^C}p^C_\Gamma\right),
\\
\Psi^{A\Gamma\Lambda}_{\alpha\beta C}:= &
\sum_{B=1}^m\left(
\frac{\d F^{A}_{\alpha} }{\d p^B_\Gamma} \frac{\d F^{B}_{\beta}}{\d p^C_\Lambda} +
\frac{\d F^{A}_{\alpha} }{\d p^B_\Lambda} \frac{\d F^{B}_{\beta}}{\d p^C_\Gamma}\right )
\end{aligned}
\end{equation}
with the notation
\begin{equation}\Label{pprime}
p=(p^A)^{A=1,\ldots,m}, \quad
p'=\left(p^A_\Gamma  \right)_{\Gamma=1,\ldots, k}^{A=1,\ldots,m}.
\end{equation}
\end{Thm}

As an illustration of Theorem \ref{main1}, we give an example where the conditions \eqref{comp10} and \eqref{comp20} take a particularly simple form.

\begin {Ex}{\rm  Let $L^A_{\alpha B}$, for $A,B=1,\ldots, m$ and $\alpha=k+1,\ldots, n$, be vector fields of the form
\begin{equation*}
L^A_{\alpha B}=\delta^A_B\frac{\d}{\d x^\alpha}+\sum_{\Lambda=1}^k\xi^{A\Lambda}_{\alpha B}(x) \frac{\d}{\d x^\Lambda},
\end{equation*}
where $\delta^A_B$ denotes the Kronecker delta and the coefficients $\xi^{A\Lambda}_{\alpha B}(x)$ are real-valued and real-analytic in an open subset $\Omega$ of $\bR^n$. Let
$\mathcal L_\alpha$ be the $m\times m$ matrix of vector fields
$\mathcal L_\alpha:=(L^A_{\alpha B})_{A,B=1,\ldots m}$ acting on $\bR^m$-valued functions $u=(u^B)_{B=1,\ldots, m}$ by
$$
\mathcal L_\alpha u:=\left (\sum_{B=1}^m L^A_{\alpha B}u^B\right )_{A=1,\ldots, m}.
$$
Clearly, the Cauchy problem for $u=(u^A)_{A=1,\ldots, m}$ given by
\begin{equation}\Label{CPF}
\left\{
\begin{aligned}
& \mathcal L_\alpha u=0\\
& u^A(x^1\ldots, x^k,x^{k+1}_0,\ldots, x^n_0)=a^A(x^1, \ldots, x^k)
\end{aligned}
\right.,
\quad A=1,\ldots, m,\ \alpha=k+1,\ldots, n,
\end{equation}
where $x^0\in\Omega$,
can be written in the form \eqref{CP0}. An elementary computation shows that the compatibility conditions \eqref{comp10} and \eqref{comp20} are equivalent to the generalized Frobenius condition
$$
[\mathcal L_\alpha,\mathcal L_\beta]:=\mathcal L_\alpha\mathcal L_\beta-\mathcal L_\beta\mathcal L_\alpha=0.
$$
Note that $\mathcal L_\alpha\mathcal L_\beta$ is an $m\times m$ matrix of second order differential operators and the commutator $[\mathcal L_\alpha,\mathcal L_\beta]$ is, in general, also an $m\times m$ matrix of second order operators. If $m=1$, then the commutator $[\mathcal L_\alpha,\mathcal L_\beta]$ is a vector field, and Theorem \ref{main1} for the Cauchy problem \eqref{CPF} reduces to the standard Frobenius theorem (in the real-analytic category).
}
\end{Ex}

As a first step in the proof of Theorem \ref{main1}, we shall prove
(i) $\implies$ (iii) (i.e.\ the necessity of the compatibility conditions)
and (ii) $\iff$ (iii) (i.e.\ the characterization of the compatibility conditions
by means of approximate solutions). This will be done in
Section \ref{solvcomp}. To complete the proof, we show,  in Section
\ref{iiitoi}, that (ii) $\implies$ (i). This step is more delicate
and uses Cartan-K\"ahler theory (see e.g.\ \cite{BCGGG}). In Sections \ref{apprsec},
\ref{ext}, and \ref{Kahler}, we introduce the necessary preliminary
material for the proof of (ii) $\implies$ (i).

The study of nonlinear overdetermined systems of analytic partial differential equations has a long history going back to the late 19th century. We mention here only some this work, namely that of E. Cartan \cite{C31}, C. Riquier  \cite{R10}, J. M. Thomas \cite{T34}, and H. Goldschmidt \cite{G67} in which results on local and formal existence of solutions to the systems are given.  The interested reader is referred to  \cite{BCGGG} for further references.

We shall apply Theorem \ref{main1} to a Cauchy problem for a system of
Monge-Amp\`ere type equations that will be used for the geometric
applications mentioned above. Let $x=(x^1,\ldots, x^n)$ denote
coordinates in $\bR^n$. For a function $u=u(x)$, we shall use the notation
$$
u_i:=u_{x^i}=\partial_{x^i}u=\frac{\partial u}{\partial x^i},
\quad
u_{ij}:=u_{x^ix^j}=\partial_{x^i}\partial_{x^j}u=\frac{\partial^2 u}{\partial x^i\partial x^j}.
$$
For integers $2\leq \alpha,\beta\leq n$, we shall denote by
$\Delta_{\alpha\beta}$ the second order, nonlinear partial
differential operator
\begin{equation}\Label{Deltaab0}
 \Delta_{\alpha\beta}(u):=\det
 \begin{pmatrix}
 u_{11}  & u_{1 \beta}\\
 u_{\alpha 1} & u_{\alpha \beta}
 \end{pmatrix}.
 \end{equation}
 We have the following result, whose proof is given in Section \ref{Monge}
 as an application of Theorem \ref{main1}.

\begin{Thm}\Label{MAthm}  Let $n>2$ and
$f_{\alpha\beta}(x,t)$, $2\leq \alpha,\beta\leq n$,
be real-analytic, real-valued  functions in a connected open subset $U\times V\subset \bR^{n}\times (\bR\setminus\{0\})$.
Then the following two conditions are equivalent:

\medskip
\noindent {\rm (i)} For any $x_0=(x^1_0,\ldots,x^n_0)\in U$ and any  real-valued, real-analytic functions $a(x^1)$, $a_2(x^1),\ldots, a_n(x^1)$ in a neighborhood of $x^1_0$ in $\bR$ such that $a''(x^1_0)\in V$, there exists a unique real-valued, real-analytic  solution $u(x)$ in a neighborhood of  $x_0$ to the
Cauchy problem
\begin{equation}\Label{MA00}
\left\{
\begin{aligned}
 & \Delta_{\alpha\beta}(u)=f_{\alpha\beta}(x,u_{11}),
 \\ &u(x^1,x^2_0,\ldots, x^n_0)=a(x^1),
 \quad u_{\alpha}(x^1,x^2_0,\ldots, x^n_0)=a_\alpha(x^1),
\end{aligned}
\right . \quad \alpha,\beta =2,\ldots, n.
\end{equation}

\medskip
\noindent
{\rm (ii)} The functions $f_{\alpha\beta}$
are of the form
$f_{\alpha\beta}(x,t)= g_{\a\b}(x)t$ with  $g_{\a\b}$
satisfying the compatibility conditions
\begin{equation}\Label{symder0}
g_{\alpha\beta}=g_{\beta\alpha},\quad \partial_{x^1}g_{\alpha\beta}=0,\quad
\partial_{x^\gamma}g_{\alpha\beta}=\partial_{x^\beta}g_{\alpha\gamma},\quad
\alpha,\beta,\gamma=2,\ldots, n.
\end{equation}
\end{Thm}

\begin{Rem}
{\rm
The conditions \eqref{symder0} are locally equivalent to
the existence of a real-analytic function $v(x)$ independent of $x_1$
with $g_{\a\b}=\d _{x_\a}\d_{x_\b} v$, for $2\le \a,\b\le n$.
Also observe that if $f_{\a\b}(x,t)$ are independent of $t$,
then (ii) is equivalent to the vanishing of all $f_{\a\b}$.}
\end{Rem}

\begin{Rem}
{\rm
Observe that when $n=2$, the condition (i) of Theorem \ref{MAthm} holds for any right hand side $f_{\a\b}=f_{22}$ as a consequence of the Cauchy-Kowalevsky Theorem. Hence, the implication (i)$\implies$(ii) does not hold for $n=2$.}
\end{Rem}

Next, recall that if $\Sigma$ is a smooth oriented hypersurface in
$\bR^{n+1}$, then its Gauss map $G\colon\Sigma \to S^n$ is the function
that sends $x\in \Sigma$ to the (positive) unit normal $n(x)\in S^n$. The
Gauss (or spherical) image of $\Sigma$ is $G(\Sigma)\subset S^n$. If
$\Sigma$ is
(locally) defined in $\bR^{n+1}$ as a graph $y=u(x)$, then the rank
of the Gauss map at a point $(x,u(x))\in \bR^n\times\bR$ is
precisely the rank of the Hessian of $u$ at $x$. Consequently, the
graphs $y=u(x)$ for which the rank of the Gauss map is identically
one (with the additional harmless assumption that $u_{11}\neq 0$)
correspond precisely to the solutions of the system of equations
$\Delta_{\alpha\beta}(u)=0$, $\alpha,\beta=2,\ldots, n$. The Cauchy
data $a(x^1)$, $a_\alpha(x^1)$, $2\leq \alpha\leq n$,  in \eqref{MA00} can be used to
prescribe the Gauss image of $\Sigma$. More precisely, we have the
following result, whose proof is given in Section \ref{realgeom} as an application of Theorem \ref{MAthm}.

\begin{Thm}\Label{Gauss}
Let $\gamma\colon (-1,1)\to S^n$ be a real-analytic curve with $\gamma(0)=(0,\ldots,0,-1)$ and $\gamma'(0)=(1,0\ldots,0)$. Then, there exists a unique (germ at $0$ of a) real-analytic function $u(x)$, $x\in\bR^n$, with $u(0)=0$
such that if $\Sigma\subset\bR^{n+1}$ is the graph of $u(x)$ and $G$ is its Gauss map, then
the rank of $G$ is equal to one at every point of $\Sigma$
and $G((x^1,0), u(x^1,0))=\gamma(x^1)$.\end{Thm}

We remark that this local result contrasts with the global
situation: Hartman and Nirenberg \cite{HN} showed that any {\it
properly embedded} smooth, connected hypersurface in $\bR^{n+1}$
whose Gauss map has rank $\leq 1$ at every point (or, equivalently,
whose Riemannian curvature tensor vanishes identically at every
point) is necessarily a cylinder over a curve, and hence its Gauss
image is contained in a 2-plane section of the sphere, i.e.\ is
contained in a circle. This strong global rigidity fails for
hypersurfaces whose Gauss map has rank two or higher. There are
relatively simple examples of properly embedded real-analytic
hypersurfaces in $\bR^{n+1}$ whose Gauss maps have rank two, but
whose Gauss images are not even contained in any hyperplane section
of the sphere (see \cite{S60} and \cite{W95}; see also \cite{AG04}).

To describe our result concerning real hypersurfaces in $\bC^{n+1}$,
we need some definitions. Let $M$ be a real-analytic, connected
hypersursurface in $\bC^{n+1}$. Recall that $M$ is said to be {\it
holomorphically nondegenerate} if there are no nontrivial germs of
holomorphic vector fields, i.e.\ holomorphic sections of
$T^{1,0}\bC^{n+1}$, that are tangent to $M$.  If $p\in M$, we let
$\rho(Z,\bar Z)$ be a local defining function for $M$ near $p$. We
let $L_1,\ldots,L_n$ be a local basis for the CR vector fields on
$M$ near $p$.  The hypersurface $M$ is said to be {\it
$k$-nondegenerate} at $p$ if the collection of vectors
\begin{equation}\Label{k-nondeg}
L^I\rho_{Z}(p,\bar p),\quad I\in \bZ_+^n,\ |I|\leq j,
\end{equation}
spans $\bC^{n+1}$ for $j=k$ and $k$ is the smallest integer with
this property; here, we use standard multi-index notation
$L^I:=L_1^{I_1}\ldots L_n^{I_n}$, etc. The notion of
$k$-nondegeneracy is independent of the choice of defining function
$\rho$ and CR vector fields $L_1,\ldots, L_n$. We shall say that $M$
is finitely nondegenerate at $p$ if it is $k$-nondegenerate at $p$
for some $k$. The reader is referred to \cite{BERbook} for further
information about these notions and CR geometry in general. We
mention that $M$ is holomorphically nondegenerate if and only if $M$
is finitely nondegenerate on a dense open subset. For a
holomorphically nondegenerate hypersurface $M$, its {\it Levi
number}, denoted $\ell(M)$, is defined to be the minimal $k$ such
that $M$ is $k$-nondegenerate at some point (and hence for all
points in a dense open subset). As mentioned above, the Levi number
was introduced in \cite{BHR}, where it was also shown that $\ell(M)$
is always $\leq n$. If $M$ is Levi nondegenerate at some point, then
$\ell(M)=1$ and, hence,  real hypersurfaces in general position satisfy
$\ell(M)=1$.
We remark that it is trivial to construct, for any
integer $k$, a real-analytic hypersurface $M$ that is
$k$-nondegenerate at a given point $p$. However, the Levi number of
such an $M$ will still be, generically, equal to one. Examples of
hypersurfaces $M$ with $\ell(M)\geq 2$ are scarce. The first example
of an $M$ in $\bC^3$ with $\ell(M)=2$ can be found in \cite{F77}.
Another example with $\ell(M)=2$ is the tube over the light cone.
 Homogeneous tube hypersurfaces  in $\bC^3$  with
Levi number $=2$ were systematically studied in \cite{FK06},
\cite{FK07}. An example in $\bC^4$ with Levi number $=3$ was
constructed in \cite{F07}. To the best of the authors' knowledge, no
examples with Levi number $\geq4$ are known in the literature. In
this paper we show that, for any $n\geq 1$, there are (plenty of)
real-analytic tube hypersurfaces in $\bC^{n+1}$ with Levi number
$=n$. This statement is an immediate consequence of Theorem
\ref{Gauss} and the following proposition; see also Remark
\ref{lastrem}.

\begin{Pro}\Label{Tube} Let $\Sigma$ be a real-analytic, connected hypersurface in
$\bR^{n+1}$. Assume that its Gauss image $G(\Sigma)$ is a real-analytic
curve in $S^n\subset \bR^{n+1}$ that is not contained in a
hyperplane.  Then, the real-analytic tube hypersurface
$M:=\Sigma+i\bR^{n+1}\subset \bC^{n+1}$ is holomorphically
nondegenerate, the rank of the Levi form of $M$ is one on a dense
open subset of $M$, and the Levi number $\ell(M)=n$.
\end{Pro}

Finally, we remark that a holomorphically nondegenerate tube
hypersurface $M:=\Sigma+i\bR^{n+1}$ in $\bC^{n+1}$ such that the rank
of its Levi form is one at every point, the existence of which is
guaranteed by Theorem \ref{Gauss} and Proposition \ref{Tube}, is
only locally defined, i.e.\ not properly embedded in $\bC^{n+1}$.
This follows again from the theorem of Hartman and Nirenberg: If $M$
were properly embedded, then $\Sigma$ would be a properly embedded,
noncylindrical hypersurface whose Gauss map has rank one at every
point, contradicting the theorem of Hartman and Nirenberg.

This paper is organized as follows. In Section \ref{CKsection}, we
introduce some basic notation and conventions that will be used
throughout the paper. The general Cauchy problem and
the main result, Theorem \ref{main1},  concerning the
existence and uniqueness of solutions are formulated. The proof of
this theorem is then given in Sections \ref{solvcomp},
\ref{apprsec}, \ref{ext}, \ref{Kahler}, and \ref{iiitoi}. The proof
of Theorem \ref{MAthm} is given in Section \ref{Monge}. The final
two sections are devoted to the proofs of Theorem \ref{Gauss} and
Proposition \ref{Tube}, respectively.

The authors wish to thank Robert Bryant for comments on an earlier version of this paper, and in particular for providing us with relevant references.

\section{Notation and conventions}\Label{CKsection}

In this section, we shall introduce some notation and conventions that will be used in this paper. Let $1\leq k<n$ be integers and
$x=(x^1,\ldots,x^k,x^{k+1},\ldots,x^n)$ denote coordinates in
$\bR^n=\bR^k\times\bR^{n-k}$. We shall use the following conventions. Small
Latin letters $i,j$, etc.\ will run over the set $\{1,\ldots,n\}$,
capital Greek letters $\Lambda,\Gamma$, etc.\ will run over
$\{1,\ldots,k\}$, small Greek letters $\alpha,\beta$, etc.\ will run
over $\{k+1,\ldots,n\}$, and capital Latin letters $A,B$, etc.\ will
run over $\{1,\ldots, m\}$. Thus, we shall refer to the coordinates
as $x^i$, and $x^\Lambda$ can be used as coordinates on the initial
$k$-plane $\{ x^\alpha=x^\alpha_0\}$ in the Cauchy problem \eqref{CP0}.  The unknown function will be denoted
$u^A$. We shall use the notation $$u^A_i:=\frac{\partial
u^A}{\partial x^i}.$$
The right hand side of the system of partial differential equations in \eqref{CP0} will be denoted by $F^A_\alpha(x^i,p^B,p^B_\Lambda)$ with the notation introduced in \eqref{pprime}. The Cauchy problem \eqref{CP0} can now be written as
\begin{equation}\Label{CP}
\left\{
\begin{aligned}
& u^A_\alpha=F^A_\alpha(x^i,u^B,u^B_\Gamma)\\
& u^A(x^\Lambda,x^\alpha_0)=a^A(x^\Lambda)
\end{aligned}
\right.,
\quad A=1,\ldots, m,\ \alpha=k+1,\ldots, n,
\end{equation}
the point $(x_0,p_0,p'_0)$ as $((x^i)_0,(p^A)_0,(p^A_\Lambda)_0)$, and the condition \eqref{initial0} as
\begin{equation}\Label{initial}
(a^A(x^\Gamma_0),a^A_\Lambda(x^\Gamma_0))=((p^A)_0,(p^A_\Lambda)_0).
\end{equation}

We shall use the following notation for the derivatives of
$F^A_\alpha$
$$
F^A_{\alpha i}:=\frac{\partial F^A_\alpha}{\partial x^i},\quad
F^A_{\alpha B}:=\frac{\partial F^A_\alpha}{\partial p^B},\quad
F^{A\Lambda}_{\alpha B}:=\frac{\partial F^A_\alpha}{\partial
p_\Lambda^B}.
$$
We will also use the summation convention that an index that appears
both as a sub- and superscript is summed over. The functions in \eqref{compfunc} appearing in the compatibility conditions \eqref{comp10} and \eqref{comp20} can now be written in the following way
\begin{equation}\Label{compfunc2}
\Phi^A_{\alpha\beta}:= F^A_{\alpha\beta}+F^A_{\alpha B}F^B_{\beta}+F^{A \Gamma}_{\alpha
B}(F^B_{\beta \Gamma}+F^B_{\beta
C}p^C_\Gamma),
\quad
\Psi^{A\Gamma\Lambda}_{\alpha\beta C}:=
F^{A \Gamma}_{\alpha B} F^{B \Lambda}_{\beta C} +
F^{A \Lambda}_{\alpha B} F^{B \Gamma}_{\beta C}.
\end{equation}

\section{Solvability and compatibility conditions}\Label{solvcomp}

We keep the notation introduced in the previous section. We begin by proving the implication (i)$\implies$(iii) in Theorem \ref{main1}.
\begin{proof}[Proof of {\rm (i)$\implies$(iii)}]
If $u^A(x^i)$ is a solution to \eqref{CP}, in particular, of the  equation
\begin{equation}\Label{PDE}
u^A_\alpha(x^i)=F^A_\alpha(x^i, u^B(x^i), u^B_{\Gamma}(x^i)),
\end{equation}
then, by differentiating with respect to $x^\Lambda$, we obtain
\begin{equation}\Label{difflambda}
u^A_{\alpha\Lambda}=F^A_{\alpha\Lambda}+F^A_{\alpha B}u^B_{\Lambda}+F^{A\Gamma}_{\alpha B}u^B_{\Gamma\Lambda}.
\end{equation}
Similarly, differentiating \eqref{PDE} with respect to $x^\beta$, we obtain
\begin{equation}\Label{diffbeta}
u^A_{\alpha\beta}=F^A_{\alpha\beta}+F^A_{\alpha B}u^B_{\beta}+F^{A\Gamma}_{\alpha B}u^B_{\Gamma\beta}.
\end{equation}
By substituting for $u^B_{\Gamma\beta}$ in \eqref{diffbeta} using \eqref{difflambda}
(and the symmetry in $\Gamma$ and $\beta$), we obtain
\begin{equation}\Label{finaldiff}
u^A_{\alpha\beta}=F^A_{\alpha\beta}+F^A_{\alpha B}u^B_{\beta}+F^{A\Gamma}_{\alpha B}(F^B_{\beta\Gamma}+F^B_{\beta C}u^C_{\Gamma}
+F^{B\Omega}_{\beta C}u^C_{\Omega\Gamma}).
\end{equation}
In particular, since the left hand side is symmetric in $\alpha$ and $\beta$, it follows that the right hand side is also symmetric. Fix a point  $((x^i)_0,(p^A)_0,(p^A_\Lambda)_0)\in U$. By assumption, the Cauchy problem \eqref{CP} has a unique solution for any Cauchy data $a^A(x^\Lambda)$ satisfying \eqref{initial}. We choose the data $a^A$ such that the second derivatives $a^A_{\Omega\Gamma}(x^\Lambda_0)$ vanish for all $A$, $\Omega$, and $\Gamma$. Evaluating \eqref{finaldiff} at $x^i=x^i_0$ and using the facts that
\begin{equation}\Label{facts}
u^A(x^i_0)=a^A(x^\L_0)=(p^A)_0, \quad
u^A_\Gamma(x^i_0)=a^A_\Gamma(x^\L_0)=(p^A_\Gamma)_0,
\quad
u^A_{\Omega\Gamma}(x^i_0)=a^A_{\Omega\Gamma}(x^\L_0)=0,
\end{equation}
we conclude that the quantity
\begin{equation}\Label{atq}
F^A_{\alpha\beta}+F^A_{\alpha B}p^B_{\beta}+F^{A\Gamma}_{\alpha B}(F^B_{\beta\Gamma}+F^B_{\beta C}p^C_{\Gamma}),
\end{equation}
evaluated at the point $((x^i)_0,(p^A)_0,(p^A_\Lambda)_0)$, is symmetric in $\alpha$ and $\beta$. Since the point $((x^i)_0,(p^A)_0,(p^A_\Lambda)_0)$ was arbitrary in $U$, we conclude that
\eqref{comp10} holds identically in $U$.

To prove \eqref{comp20}, we again fix a point  $((x^i)_0,(p^A)_0,(p^A_\Lambda)_0)\in U$ and choose the Cauchy data $a^A(x^\Lambda)$, satisfying \eqref{initial},  such that the second derivatives $a^A_{\Omega\Gamma}(x^\Lambda_0)$ are
equal to arbitrary chosen real numbers $q^A_{\Omega\Gamma}$
with $q^A_{\Omega\Gamma}$ being symmetric in $\Omega$ and $\Gamma$.
Evaluating \eqref{finaldiff} at $(x^i)_0$ and using
the symmetry of \eqref{comp10} in $\alpha$ and $\beta$ as well as \eqref{facts}, we conclude that
\begin{equation}
(F^{A\Gamma}_{\alpha B}F^{B\Omega}_{\beta C}-F^{A\Gamma}_{\beta B}F^{B\Omega}_{\alpha C})q^C_{\Omega\Gamma}=0.
\end{equation}
The symmetry condition given by  \eqref{comp20} now follows, since the
$q^C_{\Omega\Gamma}$ (symmetric in $\Omega$ and $\Gamma$) and the point
$(x^i_0,(p^A)_0,(p^A_\Lambda)_0)\in U$ are arbitrary. This completes the proof of (i)$\implies$(iii).
\end{proof}

\begin{proof} [Proof of {\rm (ii)$\iff$(iii)}]  In order to prove the equivalence (ii)$\iff$(iii), we shall need the following proposition.

\begin{Pro}\Label{eqapprox} Let $((x^i)_0,(p^A)_0,(p^A_\Lambda)_0)\in U$ and
$a^A(x^\Lambda)$ be real-analytic initial data defined in a neighborhood of $(x^\Lambda)_0$
satisfying \eqref{initial}.
Then $u^A(x^i)$ is an approximate solution to \eqref{CP} at $(x^i)_0$ if and only if the following equalities hold:
\begin{equation}\Label{altappsol}
\begin{aligned}
u^A &=a^A,\quad
u^A_\Gamma=a^A_\Gamma,\quad
u^A_{\Gamma\Omega}=a^A_{\Gamma\Omega},\\
u^A_\alpha &= F^A_\alpha\\
u^A_{\alpha\Lambda} &= F^A_{\alpha\Lambda}+
F^A_{\alpha B} a^B_{\Lambda}+F^{A\Gamma}_{\alpha B} a^B_{\Gamma\Lambda}
\\
u^A_{\alpha\beta} &= F^A_{\alpha\beta}+F^A_{\alpha B}F^B_{\beta}+F^{A \Gamma}_{\alpha
B}(F^B_{\beta \Gamma}+
F^B_{\beta
C}a^C_\Gamma +F^{B\Omega}_{\beta C}a^{C}_{\Omega\Lambda}),
\end{aligned}
\end{equation}
where all the $u^A$ (resp.\ $a^A$, resp.\ $F^A_\alpha$) and their derivatives appearing in the right hand sides of \eqref{altappsol} are evaluated at the point $(x^i)_0$ (resp.\ $(x^\L)_0$, resp.\ $((x^i)_0,(p^A)_0,(p^A_\Lambda)_0)$).
\end{Pro}

\begin{proof}  Clearly, the first row of equations in \eqref{altappsol} is equivalent to the initial condition in \eqref{CPapprox}. Now, proceeding as in the proof of (i)$\implies$ (iii) above, it is straightforward to verify that $u^A(x^i)$ satisfies the  first row in \eqref{CPapprox} if and only if the remaining equations in \eqref{altappsol} hold. This completes the proof of Proposition \ref{eqapprox}.
\end{proof}

The following proposition follows immediately from Proposition \ref{eqapprox} by similar arguments to those used to prove  (i)$\implies$ (iii) above.

\begin{Pro}\Label{approxtocomp} Let $((x^i)_0,(p^A)_0,(p^A_\Lambda)_0)\in U$. Then, the following are equivalent:
\medskip

\noindent {\rm (ii$'$)} The Cauchy problem \eqref{CP} is approximately
solvable for every initial data at $((x^i)_0,(p^A)_0,(p^A_\Lambda)_0)$.

\medskip
\noindent {\rm (iii$'$)} The functions $F^A_\alpha$ satisfy the
compatibility conditions \eqref{comp10} and \eqref{comp20} at the point $((x^i)_0,(p^A)_0,(p^A_\Lambda)_0)$.

\end{Pro}

Clearly, the equivalence (ii)$\iff$(iii) in Theorem~\ref{main1} is a direct consequence of Proposition \ref{approxtocomp} by letting the point  $((x^i)_0,(p^A)_0,(p^A_\Lambda)_0)$ vary in $U$.
\end{proof}

\section{Approximate solvability with respect to a slope}\Label{apprsec}

We begin this section by considering a somewhat more general system of first order partial differential equations, using still the conventions introduced in previous sections:
\begin{equation}\Label{generalsystem}
G^A_\alpha(x^i,u^B,u^B_i)=0, \quad \alpha=k+1,\ldots, n,\ A=1,\ldots, m,
\end{equation}
where $G^A_\alpha$ are real-analytic functions defined in an open subset $V$ of $\bR^n\times\bR^m\times \bR^{nm}$ and $((x^i)_0,(p^A)_0,(p^A_i)_0)\in V$ is such that
$G^A_\alpha((x^i)_0,(p^A)_0,(p^A_i)_0)=0$.  Let $M$ a $k$-dimensional (real-analytic) submanifold in $\bR^n$ through $x_0$ whose tangent $k$-plane $P$ at $x_0$ is given by
\begin{equation}\Label{plane}
P=T_{x_0}M=\{x^i\colon (x^j-x^j_0)\xi^\beta_j=0,\ \beta=k+1,\ldots,n\},
\end{equation}
where $\xi^\beta:=(\xi^\beta_i)_{i=1,\ldots,n}$ are linearly independent vectors in $\bR^n$.
We shall say that $M$ is {\it non-characteristic} at $((x^i)_0,(p^A)_0,(p^A_i)_0)$ for the system \eqref{generalsystem} if the $m(n-k)\times m(n-k)$ matrix $W:=(W^{A,\beta}_{B,\alpha})$, where $A,B$ run over $\{1,\ldots m\}$ and $\alpha,\beta$ over $\{k+1,\ldots, n\}$, is invertible with
\begin{equation}\Label{W}
W^{A,\beta}_{B,\alpha}:=G^{A j}_{\alpha B}((x^i)_0,(p^A)_0,(p^A_i)_0)\xi^\beta_j,
\end{equation}
where we continue using our convention that adding extra indices corresponds to derivatives,
i.e.\ $G^{A j}_{\alpha B}$ is the derivative of $G^{A }_{\alpha}$ in $p^B_j$.

\begin{Pro}\Label{charac} Suppose that a $k$-dimensional real-analytic submanifold $M$ of $\bR^n$
through $x_0$ is non-characteristic at $((x^i)_0,(p^A)_0,(p^A_i)_0)$ for the system \eqref{generalsystem}. Then after a real-analytic change of coordinates in $\bR^n$
near $x_0$ preserving $x_0$,
$M$ can be written near $x_0$ as $\{x\colon x^\alpha=x^\alpha_0,\ \alpha=k+1,\ldots, n\}$ and
the system \eqref{generalsystem} near $((x^i)_0,(p^A)_0,(p^A_i)_0)$ as
\begin{equation}\Label{Fsystem}
u^A_\alpha=F^A_\alpha(x^i,u^B,u^B_\Lambda),\quad A=1\ldots, m,\ \alpha=k+1,\ldots n,
\end{equation}
where $F^A_\alpha(x^i,p^B,p^B_\Lambda)$ is real-analytic near the point
$((x^i)_0,(p^A)_0,(p^A_\L)_0)$ with
$(p^A_\alpha)_0=F^A_\alpha((x^i)_0,(p^B)_0,(p^B_\Lambda)_0)$.
\end{Pro}

\begin{proof} Let $(\xi^\b)$ be as in \eqref{plane} and $(\xi^\L)$ be further vectors in $\bR^n$ such that the collection $(\xi^i)$ forms a basis for $\bR^n$ and let $\xi$ denote the invertible $n\times n$ matrix whose $i$-th row equals $\xi^i$.
Then there exists a real-analytic change of coordinates near $x_0$ of the form
\begin{equation}
\tilde x^i:=x^i_0+\xi^i_j(x^j-x^j_0)+O(\sum_j|x^j-x^j_0|^2),
\end{equation}
such that $M$ near $x_0$ is given by $\tilde x^\a=\tilde x^\a_0$, where $\tilde x^i_0=x^i_0$. Let $\eta$ denote the inverse of the matrix $\xi$. Then, an application of the chain rule shows that in the new coordinates, the system \eqref{generalsystem} takes the form $\tilde G^A_\alpha(\tilde x^i,\tilde u^B,\tilde u^B_i)=0$, where
\begin{equation}
\tilde G^A_\alpha(\tilde x^i,\tilde p^B,\tilde p^B_i)=
G^A_\alpha(x^i_0+\eta^i_j(\tilde x^j-x^j_0), \tilde p^B, \xi^j_i\tilde p^B_j)+O(\sum_j|x^j-x^j_0|^2).
\end{equation}
Note that $\tilde G^{A\beta}_{\alpha B}((\tilde x^i)_0,(\tilde p^A)_0,(\tilde p^A_i)_0)=W^{A,\beta}_{B,\alpha}$, the latter  given by \eqref{W}. Thus, the conclusion of Proposition \ref{charac} follows immediately from the implicit function theorem.
\end{proof}

To prove  (ii)$\implies$(i) in Theorem \ref{main1}, we shall need some further preliminary material.  Let $F^A_\alpha$ be as in previous sections,  $((x^i)_0,(p^A)_0,(p^A_\Lambda)_0)$ be a point in $U$, $(c^\alpha_\Lambda)$ be a real $(n-k)\times k$ matrix, and $P\subset \bR^n$ the $k$-plane through the point $x_0$ given by
\begin{equation}\Label{planeP}
P:=\{x : x^\alpha=x^\alpha_0+c^\alpha_\Lambda(x^\Lambda-x^\Lambda_0),  \,\alpha=k+1,\ldots, n\}.
\end{equation}
Note that $P$, given by \eqref{planeP}, is non-characteristic at
$((x^i)_0,(p^A)_0,(p^A_i)_0)$ with $(p^A_\a)_0=F^A_\a((x^i)_0,(p^A)_0,(p^A_\L)_0)$
for the system \eqref{Fsystem} if the $m(n-k)\times m(n-k)$ matrix $V:=(V^{A,\beta}_{B,\alpha})$ (where as before $A,B$ run over $\{1,\ldots m\}$ and $\alpha,\beta$ over $\{k+1,\ldots, n\}$) is invertible with
\begin{equation}\Label{V}
V^{A,\beta}_{B,\alpha}:=\delta^A_B\delta^\beta_\alpha-F^{A\Lambda}_{\alpha B}
((x^i)_0,(p^A)_0,(p^A_\Lambda)_0)c^\beta_\Lambda,
\end{equation}
where $\delta^i_j$ is the standard Kronecker delta symbol. In this case, we shall say that the slope $(c^\alpha_\Lambda)$ is non-characteristic at $((x^i)_0,(p^A)_0,(p^A_\Lambda)_0)$ for the system \eqref{Fsystem}.

For non-characteristic slopes $(c^\alpha_\Lambda)$, we shall consider the Cauchy problem for
\begin{equation}\Label{eqn}
u_\alpha ^A=F^A_\alpha(x^i,u^B, u^B_\Lambda)
\end{equation}
 with data on the tilted $k$-plane \eqref{planeP}.
We shall say that this Cauchy problem with respect to the (non-characteristic) slope $(c^\alpha_\Lambda)$ is {\it approximately solvable for every initial data at $((x^i)_0,(p^A)_0,(p^A_\L)_0)$} if, for every real-analytic $a^A(x^\Lambda)$ defined in a neighborhood
of $(x^\L)_0$ and satisfying
 \begin{equation}
\Label{firstjet}
a^A(x^\Gamma_0)=(p^A)_0, \quad
a^A_\Lambda(x^\Gamma_0)=(p^A_\Lambda)_0+
F^A_\alpha((x^i)_0,(p^A)_0,(p^A_\Omega)_0)c^\alpha_\Lambda,
\end{equation}
there is an approximate solution $u^A(x^i)$ at $x_0$ to the Cauchy problem
\begin{equation}\Label{CPjettilt}
\left\{
\begin{aligned}
& u^A_\alpha=F^A_\alpha(x^i,u^B,u^B_\Gamma)
\\
& u^A(x^\Lambda,x^\beta_0+c^\beta_{\Gamma}(x^\Gamma-x^\Gamma_0))=a^A(x^\Lambda).
\end{aligned}
\right.
\end{equation}
(That is, the first equation in \eqref{CPjettilt} holds up to order $1$ at $(x^i)_0$ and the second up to order $2$ at $(x^\L)_0$.)
Note that \eqref{firstjet} is a modification of \eqref{initial}
corresponding to the tilted plane \eqref{planeP}.
This choice also has the effect that the definition of approximate solvability with respect to a slope is invariant under linear changes of coordinates.
In the terminology of the previous section, the approximate solvability of the Cauchy problem
corresponds to that with respect to the slope $(c^\alpha_\Lambda)$ with $c^\alpha_\Lambda=0$.

The following statement relates approximate solvability with respect to different slopes
and is of independent interest.

\begin{Pro}\Label{slopeprop} The Cauchy problem for \eqref{eqn} with respect to a fixed non-characteristic slope is approximately solvable for every initial data at
$((x^i)_0,(p^A)_0,(p^A_\L)_0)$  if and if it is approximately solvable for every initial data at $((x^i)_0,(p^A)_0,(p^A_\L)_0)$ with respect to every non-characteristic slope.
\end{Pro}

\begin{proof} We observe that all results proved so far also hold in the complex-analytic setting, i.e.\ by letting the variables $(x^i,p^B, p^B_\Lambda)$ be complex and the functions $F^A_\alpha$, $u^A$, and the initial data $a^A$ be holomorphic with respect to their arguments. In this proof, we shall first consider this complex setting. Thus, all variables will be understood to be complex and all functions holomorphic. Let $a^A(x^\Lambda)$ be any data satisfying \eqref{firstjet}.  As in Proposition \ref{eqapprox} and its proof, a function  $u^A(x^i)$ is an approximate solution to \eqref{CPjettilt} at $((x^i)_0,(p^A)_0,(p^A_\L)_0)$ (for $a^A$ satisfying \eqref{firstjet}) if and only if the following equalities hold:
\begin{equation}\Label{altappsolc}
\begin{aligned}
u^A &=a^A,\quad
u^A_\Gamma + F^A_\a c^\a_\Gamma=a^A_\Gamma,\quad
u^A_{\Lambda\Gamma}+u^A_{\Lambda\b}c^\b_{\Gamma}+u^A_{\alpha\Gamma}c^\alpha_\Lambda+u^A_{\alpha\beta}c^\alpha_{\Lambda}c^{\beta}_{\Gamma}=a^A_{\Lambda\Gamma}
\\
u^A_\alpha &= F^A_\alpha\\
u^A_{\alpha\Lambda} &= F^A_{\alpha\Lambda}+F^A_{\alpha B} (p^B_{\Lambda})_0+F^{A\Gamma}_{\alpha B} u^B_{\Gamma\Lambda}
\\
u^A_{\alpha\beta} &= F^A_{\alpha\beta}+F^A_{\alpha B}F^B_{\beta}+F^{A \Gamma}_{\alpha
B}(F^B_{\beta \Gamma}+F^B_{\beta
C}(p^C_\Gamma)_0 +F^{B\Omega}_{\beta C}u^{C}_{\Omega\Lambda}),
\end{aligned}
\end{equation}
where all the $u^A$ (resp.\ $a^A$, resp.\ $F^A_\alpha$) and their derivatives appearing in the right hand sides are evaluated at the point $(x^i)_0$ (resp.\ $(x^\L)_0$, resp.\ $((x^i)_0,(p^A)_0,(p^A_\Lambda)_0)$).
We observe that, for small values of the slope $c^\alpha_{\Lambda}$, there is always a unique solution to \eqref{altappsolc} consisting of the value of $u^A$ and its derivatives up to order $2$ at $x_0$. To see this, notice that the system is trivial to solve when $c^\alpha_\Lambda=0$.
As in the proof of (i)$\implies$(iii) above (cf.\ Proposition \ref{approxtocomp}),
we see that for $(c^\alpha_\Lambda)$ sufficiently small, \eqref{altappsolc}
implies \eqref{comp10} and \eqref {comp20}. Indeed $u_{\a\b}=u_{\b\a}$
implies the symmetry in $\a$ and $\b$ of the right-hand side of the last identity
in \eqref{altappsolc}. Furthermore, for $c^\a_\L$ small, as the $a^A_{\Lambda\Gamma}$ vary over all possible choices, so do the $u^A_{\Lambda\Gamma}$ satisfying \eqref{altappsolc}.
Hence we obtain the compatibility conditions \eqref{comp10} and \eqref {comp20} at $((x^i)_0,(p^A)_0,(p^A_\Lambda)_0)$. Vice versa, if \eqref{comp10} and \eqref {comp20} hold at $((x^i)_0,(p^A)_0,(p^A_\Lambda)_0)$, the right-hand side of the last identity
in \eqref{altappsolc} is always symmetric in $\a$ and $\b$.
Hence the system \eqref{altappsolc} is always solvable in $u^A_{\L\Gamma}$, $u^A_{\a\L}$, $u^A_{\a\b}$ for $c^\a_\L$ small.

In particular, we conclude that approximate solvability for every initial data at
$((x^i)_0,(p^A)_0,(p^A_\Lambda)_0)$ with respect to one sufficiently small
slope $(c^\alpha_\Lambda)$ is equivalent to approximate solvability for every initial data at
$((x^i)_0,(p^A)_0,(p^A_\Lambda)_0)$ with respect to every sufficiently small slope $(c^\alpha_\Lambda)$.
We claim that this fact implies that the set $\mathcal S$ of all non-characteristic complex slopes $(c^\alpha_\Lambda)$ with respect to which approximate solvability (for every initial data) holds at $((x^i)_0,(p^A)_0,(p^A_\Lambda)_0)$ is open and closed in the space of all non-characteristic slopes. Indeed, given any non-characteristic $(c^\alpha_\Lambda)$, then, after a linear change of coordinates of the form $\tilde x^\Lambda=x^\Lambda$,
$\tilde x^\alpha=x^\alpha-c^\alpha_\Lambda(x^\Lambda-x^\Lambda_0)$, the initial $k$-plane \eqref{planeP} becomes $\{\tilde x^\alpha=x^\alpha_0\}$ and the system of differential equation
\eqref{eqn} can be written as $\tilde u^A_\alpha=\tilde F^A_\alpha(\tilde x^i,\tilde u^B, \tilde u^B_\Lambda)$ by Proposition \ref{charac}.
Since approximate solvability with respect to slopes is invariant under affine changes of coordinates,
the observation above shows that $\mathcal S$ is both open and closed in the space of all non-characteristic slopes. Now, observe that the set of all non-characteristic slopes is Zariski-open in the space of all complex $(n-k)\times k$ matrices. In particular, this set is connected, which by the above proves that the set $\mathcal S$ is either empty or equals the whole set of non-characteristic complex slopes. By specializing to real non-characteristic slopes, the proposition follows.
\end{proof}

\section{Exterior differential systems and integral manifolds}\Label{ext}

In this section, we shall reformulate solvability of the Cauchy
problem \eqref{CP} in terms of the existence and uniqueness of
integral manifolds for an exterior differential system. The reader
is referred to the text \cite{BCGGG} for further details of the theory of exterior
differential systems and their integral manifolds.
We here briefly recall the needed terminology in the real-analytic case.
On a real-analytic manifold $X$,
consider the graded ring of all real-analytic exterior differential forms
(of any degree) on $X$ with respect to the exterior product.
That is, an element of the ring is a (formal) finite sum of exterior forms of different degrees.
A real-analytic {\em exterior differential system}
 is any differential ideal  $\mathcal I$ in this ring,
 i.e.\ any algebraic ideal
that is closed under exterior differentiation.
A real submanifold $S\subset X$ is an {\em integral submanifold} of  $\mathcal I$
if every element from  $\mathcal I$ vanishes when restricted to $S$.
Here a formal sum of exterior differential forms of different degrees is said to vanish on $S$
if each homogeneous component of the formal sum does.

Now consider a solution  $u^A(x^i)$ to the system of
partial differential equations
\begin{equation}\Label{PDE2}
u^A_\alpha(x^i)=F^A_\alpha(x^i, u^B(x^i), u^B_{\Gamma}(x^i)),
\end{equation}
with $F^A_\a$ being real-analytic in an open subset $U\subset \bR^n\times\bR^m\times\bR^{km}$.
Then its jet-graph
\begin{equation}\Label{jetgraph}
S_{u^A}:=\{(x^i,p^B,p^B_j)\colon p^B=u^B(x^i),\ p^B_j=u^B_j(x^i)\}
\subset \bR^n\times\bR^m\times\bR^{nm},
\end{equation}
is an integral submanifold for the exterior differential system $\mathcal I$ generated by the forms
\begin{equation}\Label{E}
p^A_\alpha-F^A_\alpha(x^i, p^B, p^B_{\Gamma}), \quad dp^B-p^B_i dx^i,
\end{equation}
defined in
\begin{equation}\Label{Utilde} \tilde U:=U\times \bR^{(n-k)m}\subset \bR^n\times\bR^m\times\bR^{nm}.
\end{equation}
Conversely, if $S\subset \tilde U$ is an integral submanifold of $\mathcal I$ such that $dx^1\wedge \ldots\wedge dx^n\neq 0$ everywhere on $S$, then $S$ is automatically $n$-dimensional and is locally a jet-graph $S_{u^A}$ of a solution $u^A(x^i)$ to \eqref{PDE2}.

The differential ideal $\mathcal I$ is generated algebraically by the differential forms
\begin{equation}\Label{EDS}
\left\{\begin{aligned}
f^A_\alpha:= & p^A_\alpha-F^A_\alpha(x^i, p^B, p^B_{\Gamma}),\\
\omega^A_\alpha:= & df^A_\alpha,\\
\eta^A:= & dp^A-p^A_idx^i,\\
\Omega^A:= & -d\eta^A= dp^A_i\wedge dx^i.
\end{aligned}
\right.
\end{equation}
For convenience of notation, we shall write $z=(x^i,p^B,p^B_i)$ and hence $z_0$ will denote the point $((x^i)_0,p^B_0,(p^B_i)_0)$.

Recall that if $E$ is a linear $l$-dimensional subspace of $T_z\bR^{n+m+nm}$, then $E$ is called an {\it integral element} of $\mathcal I$ if for all $p$-forms $\phi\in\mathcal I$, we have $\phi(v_1,\ldots,v_p)=0$ for all $v_1,\ldots, v_p\in E$. If $E$ is an integral element, then its {\it polar space} (with respect to the differential ideal $\mathcal I$), denoted by $H(E)$, is given by
\begin{equation}\Label{polar0}
H(E):=\{v\in T_z\bR^{n+m+nm}\colon \phi(v,e_1,\ldots,e_l)=0,\ \forall \text{{\rm $(l+1)$-forms $\phi$ in $\mathcal I$}}\},
\end{equation}
where $e_1,\ldots, e_l$ is a basis for $E$.

\begin{Rem}\Label{int-rem}
It follows directly from the definition that if $E$ is an integral element, it is always contained in $H(E')$ for any $E'\subset E$.
\end{Rem}
We shall use the following elementary lemma, whose proof is left to the reader.

\begin{Lem}\Label{polar}
Let $E\subset T_z\bR^{n+m+nm}$ be an integral element of $\mathcal I$.
Then, a vector $v$ is in  $H(E)$ if and only if
\begin{equation}\Label{polarsystem}
\omega^A_\alpha(v)=0,\ \eta^A(v)=0,\ \Omega^A(v,e)=0, \quad \forall e\in E.
\end{equation}
\end{Lem}

We shall also need the following proposition.

\begin{Pro} \Label{HEn} Let $E\subset T_z\bR^{n+m+nm}$ be an integral
element. Assume that there are vectors $v_1,\ldots,v_k\in E$ such
that
\begin{equation}\Label{wedge}
(dx^1\wedge\ldots\wedge dx^k)(v_1,\ldots,v_k)\neq 0
\end{equation}
and such that the push forwards of these vectors by the projection
$\pi \colon \bR^{n+m+nm}\to \bR^n$, given by
$\pi(x^i,p^A,p^A_j)=(x^i)$, span a non-characteristic $k$-plane of \eqref{PDE2} at $z$. Then $\dim H(E)\leq n$.
\end{Pro}

\begin{proof} The assumption \eqref{wedge} is equivalent to the fact that there are $k$ vectors in the span of $v_1,\ldots, v_k$ of the form
\begin{equation}\Label{v1tovk}
e_\Lambda=\frac{\partial}{\partial x^\Lambda}+e^{\alpha}_\Lambda
\frac{\partial}{\partial x^\alpha}+
e_\Lambda^A\frac{\partial}{\partial p^A} + e_{\Lambda
j}^A\frac{\partial} {\partial p_j^A}.
\end{equation}
The assumption that the span of $v_1,\ldots,v_k$ projects onto a
noncharacteristic $k$-plane means that the slope
$(c^\alpha_\Lambda):=(e^\alpha_\Lambda)$ is noncharacteristic at
$z$, i.e.\ the matrix \eqref{V} is invertible with
$(c^\alpha_\Lambda):=(e^\alpha_\Lambda)$. Let $v$ be given by
\begin{equation}
v=d^j\frac{\partial}{\partial x^j}+d^A\frac{\partial}{\partial p^A} + d_{ j}^A\frac{\partial}{\partial p_j^A}.
\end{equation}
By Lemma \ref{polar}, if $v$ belongs to $H(E)$, then \eqref{polarsystem} holds, in particular, for all $e=e_\Lambda$. A straightforward calculation shows that the equations \eqref{polarsystem} with $e=e_\Lambda$ are given by
\begin{equation}\Label{HEeqs}
\begin{aligned}
\omega^A_\alpha(v) &= d^A_\alpha-d^jF^A_{\alpha j} -d^B F^A_{\alpha B} -d^B_\Lambda F^{A \Lambda}_{\alpha B}=0\\ \eta^A(v) &=d^A-d^ip^A_i=0\\
\Omega^A(v,e_\Lambda) &=d^A_\Lambda + d^A_\alpha e^\alpha_\Lambda - e^A_{\Lambda i}d^i=0.
\end{aligned}
\end{equation}
If we now solve for $d^A$ in the second line of \eqref{HEeqs}, for
$d^A_\Lambda$ in the third line, and substitute the result in the
first line, we obtain equations for $d^A_\alpha$ of the form
\begin{equation}\Label{HEeq1}
(\delta_B^A\delta_\alpha^\beta-F^{A\Lambda}_{\alpha
B}e^\beta_\Lambda)d^B_\beta=R^A_\alpha,
\end{equation}
where the $R^A_\alpha$ depend only on the $d^i$. Since the matrix
\eqref{V} is invertible with
$(c^\alpha_\Lambda):=(e^\alpha_\Lambda)$, we conclude that we may
solve \eqref{HEeq1} for $d^A_\alpha$ in terms of the $d^i$. By
substituting back into the previously solved equations, we conclude
that the coefficients $d^A$, $d^A_i$ are all determined by the
coefficients $d^i$. This proves that $\dim H(E)\leq n$.
\end{proof}

For $z_0:=((x^i)_0,p^B_0,(p^B_i)_0)\in \tilde U$ (given by \eqref{Utilde}), let $u^A(x)$ be a real-analytic function satisfying $u^A(x_0)=(p^A)_0$ and $u^A_j(x_0)=(p^A_j)_0$.
Denote again by $S_{u^A}$ its jet-graph, given by \eqref{jetgraph}. Observe that the exterior forms $\eta^A$ and $\Omega^A=d\eta^A$ in \eqref{EDS} vanish identically on the (any) jet-graph $S_{u^A}$. Also, note that $u^A$ is an approximate solution to the system \eqref{PDE2} at $x_0$ if and only if the restriction of the functions $f^A_\alpha$ in \eqref{EDS} to $S_{u^A}$ vanishes up to second order at  $z_0$ (see \eqref{CPapprox}). Consequently, $u^A$ is an approximate solution to the system \eqref{PDE2} at $x_0$ if and only if the tangent space $T_{z_0}S_{u^A}$ is an integral element of the differential ideal $\mathcal I$ generated by \eqref{EDS}.
Summarizing and using Remark~\ref{int-rem}, we obtain the following proposition.

\begin{Pro} \Label{HES} Let $z_0:=((x^i)_0,p^B_0,(p^B_i)_0)\in \tilde U$ and assume that $S_{u^A}$ is a jet-graph through $z_0$ of an approximate solution $u^A$ to \eqref{PDE2} at $x_0$.
Then $T_{z_0}S_{u^A}$ is an integral element of $\mathcal I$.
If $E\subset T_{z_0}S_{u^A}$, then $T_{z_0}S_{u^A}\subset H(E)$.
\end{Pro}

Combining Propositions \ref{HEn} and \ref{HES}, we obtain the following:

\begin{Cor}\Label{HE} Let $z_0:=((x^i)_0,p^B_0,(p^B_i)_0)\in \tilde U$ and assume that $S_{u^A}$ is a jet-graph through $z_0$ of an approximate solution $u^A$ to \eqref{PDE2} at $x_0$. If $E\subset T_{z_0}S_{u^A}$ and the assumption of Proposition $\ref{HEn}$ holds, then $T_{z_0}S_{u^A}= H(E)$.
\end{Cor}

Before entering the proof of {\rm (ii)}$\implies$  {\rm (i)} in Theorem \ref{main1},
we shall need one more result.

\begin{Pro}\Label{EinTS} Assume that property {\rm (ii)} of Theorem $\ref{main1}$
holds. Let $z_0:=((x^i)_0,p^B_0,(p^B_i)_0)\in \tilde U$ and
$E\subset T_{z_0}\bR^{n+m+nm}$ be an integral element of $\mathcal I$. Assume that there are
vectors $v_1,\ldots, v_k\in E$ such that \eqref{wedge} holds and such that the push
forwards of these vectors by the projection $\pi \colon \bR^{n+m+nm}\to \bR^n$, given
by $\pi(x^i,p^A,p^A_j)=(x^i)$, span a non-characteristic $k$-plane of \eqref{PDE2} at $z_0$. Then
there exists an approximate solution $u^A(x^i)$ to \eqref{PDE2} at
$x_0$ such that its jet graph $S_{u^A}$ contains $z_0$ and $E\subset T_{z_0} S_{u^A}$.
\end{Pro}

\begin{proof} After replacing the $k$ vectors $v_1,\ldots, v_k$ by suitable linear combinations, we may assume that there are vectors $e_1,\ldots, e_k$ of the form \eqref{v1tovk} whose push forward via $\pi$ span a non-characteristic $k$-plane in $\bR^n$. The latter condition is equivalent to the slope $(c^\alpha_\Lambda)$, given by $c^\alpha_\Lambda=e^\alpha_\Lambda$, being non-characteristic at $z_0$. As in the proof of Proposition \ref{slopeprop}, we make the linear change of coordinates $\tilde x^\Lambda=x^\Lambda$, $\tilde x^\alpha=x^\alpha-c^\alpha_\Lambda(x^\Lambda-x^\Lambda_0)$ and observe that, in the new coordinates, the system of differential equations \eqref{PDE2} can be written as $\tilde u^A_\alpha=\tilde F^A_\alpha(\tilde x^i,\tilde u^B, \tilde u^B_\Lambda)$ by Proposition \ref{charac} and the Cauchy problem for this equation is still approximately solvable for every initial data at $((x^i)_0,(p^B)_0,(p^B_\Lambda)_0)$ in view of Proposition \ref{slopeprop}. Hence, without loss of generality, we may assume that the $k$ vectors $e_\Lambda$ are of the form
\begin{equation}\Label{el}
e_\Lambda=\frac{\partial}{\partial x^\Lambda}+
e_\Lambda^A\frac{\partial}{\partial p^A} + e_{\Lambda j}^A\frac{\partial}{\partial p_j^A}.
\end{equation}
The fact that $E$ is an integral element implies, using the notation \eqref{EDS}, that we have the relations
\begin{equation}\Label{elambda}
\begin{aligned}
\eta^A(e_\Lambda) &=e^A_\Lambda-(p^A_\Lambda)_0=0\\
\Omega^A(e_\Lambda,e_\Gamma) &= e^A_{\Lambda\Gamma}-e^A_{\Gamma\Lambda}=0\\
\omega^A_\alpha(e_\Lambda)&= e^A_{\Lambda\alpha}-F^A_{\alpha\Lambda}-F^A_{\alpha B}e^B_\Lambda-F^{A\Gamma}_{\alpha B}e^B_{\Lambda\Gamma}=0.
\end{aligned}
\end{equation}
We now consider the Cauchy problem \eqref{CP} with initial data $a^A$ satisfying \eqref{initial} and $a^A_{\Lambda\Gamma}((x^\L)_0)=e^A_{\Lambda\Gamma}$.
By assumption, this Cauchy problem has an approximate solution $u^A(x^i)$.
Hence the identities \eqref{altappsol} hold (when evaluated at suitable points as outlined below \eqref{altappsol}). Let $S_0$ and $S_{u^A}$ denote the corresponding jet graphs, i.e.\ the two manifolds parametrized near $z_0$, respectively,  by the two maps $\Phi\colon \bR^k\to \bR^{n+m+nm}$ and $\Psi\colon \bR^n\to \bR^{n+m+nm}$, given by
\begin{equation}\Label{Phi}
\Phi(x^\Lambda):=\big((x^\Lambda,x^\alpha_0),a^A(x^\Lambda),(a^A_\Gamma(x^\Lambda), F^A_\alpha((x^\Lambda,x^\alpha_0),a^B(x^\Lambda),a^B_\Gamma(x^\Lambda)))\big),
\end{equation}
\begin{equation}\Label{Psi}
\Psi(x^i):=(x^i,u^A(x^i),u^A_j(x^i)).
\end{equation}
In the following we consider standard tangent vectors $\d/\d x^i\in T_{x_0}\bR^n$.
Making use of the equations \eqref{elambda}, it is straightforward to check that we have $$\Phi_*(\partial/\partial x^\Lambda)=\Psi_*(\partial/\partial x^\Lambda)=e_\Lambda,\quad \Lambda=1,\ldots, k,$$
with  $e_\L$ being given by \eqref{el}.
Furthermore, we have
\begin{equation}\Label{ealpha}
e_\alpha:=\Psi_*(\partial/\partial x^\alpha)=
\frac{\partial }{\partial x^\alpha}+u^A_\alpha
\frac{\partial }{\partial p^A}+u^A_{\alpha\Lambda}\frac{\partial }{\partial p^A_\Lambda}
+u^A_{\alpha\beta}\frac{\partial }{\partial p^A_\beta},
\end{equation}
where the coefficients $u^A_\alpha$, $u^A_{\alpha j}$ are evaluated at $x_0$ and satisfy the equations in \eqref{altappsol}. By the construction, the $n$ vectors $e_1,\ldots, e_n$ from \eqref{el} and \eqref{ealpha} form a basis for $T_{z_0}S_{u^A}$. In order to show that $E\subset T_{z_0}S_{u^A}$, we let $w$ be a vector in $E$, which is not contained in the span of $e_1,\ldots, e_k$. Without loss of generality, we may assume that $w$ is of the form
\begin{equation}\Label{w0}
w=w^\alpha\frac{\partial }{\partial x^\alpha}+w^A
\frac{\partial }{\partial p^A}+w^A_{\Lambda}\frac{\partial }{\partial p^A_\Lambda}+w^A_{\alpha}\frac{\partial }{\partial p^A_\alpha},
\end{equation}
for some real coefficients $w^\alpha$, $w^A$, $w^A_j$.
The fact that $E$ is an integral element implies, in view of \eqref{EDS}, that
\begin{equation}\Label{w}
\begin{aligned}
\eta^A(w) &=w^A-(p^A_\beta)_0 w^\beta=0,\\
\Omega^A(w, e_\Lambda) &= w^A_\Lambda-e^A_{\Lambda\beta}w^\beta=0,\\
\omega^A_\alpha(w)&= w^A_{\alpha}-F^A_{\alpha\beta}w^\beta-F^A_{\alpha B}w^B-F^{A\Gamma}_{\alpha B}w^B_{\Gamma}=0.
\end{aligned}
\end{equation}
In Section~\ref{solvcomp} we have proved that (ii) implies (iii) in Theorem~\ref{main1}, i.e.\ we may assume that the compatibility conditions \eqref{comp10} and \eqref{comp20} hold.
Now, a straightforward but tedious computation, using \eqref{altappsol} to substitute in \eqref{ealpha}, and \eqref{elambda} and \eqref{w} to substitute in \eqref{w0}, shows that $w=w^\alpha e_\alpha$, which proves the desired inclusion $E\subset T_{z_0}S_{u^A}$.
\end{proof}

\section{K\"ahler ordinary and K\"ahler regular integral
elements}\Label{Kahler}

In this section, we recall the notions of K\"ahler ordinary and K\"ahler regular integral elements for an exterior differential system (see \cite{BCGGG}) and verify them for the ideal $\mathcal I$ introduced in the previous section. Let $t_0\in \bR^N$ and $\mathcal F$ a family of germs at $t_0$ of real-analytic functions vanishing at $t_0$. Recall that $t_0$ is called an ordinary zero of $\mathcal F$ if there are $f_1,\ldots f_\kappa\in\mathcal F$ such that
\begin{equation}
\Rk \left(\frac{\partial f_i}{\partial t^j}(t_0)\right)=\kappa,
\end{equation}
where $i=1,\ldots , \kappa$, and $j=1,\ldots,N$,
and the germ at $t_0$ of the real-analytic set $\{t\colon f(t)=0,\ \forall f\in \mathcal F\}$ coincides with the germ at $t_0$ of  $\{t\colon f_1(t)=\ldots =f_\kappa(t)=0\}$.

 For a fixed $l\geq k$, we shall use the following convention. Lower-case $a,a',a'',\ldots$ run over the set of indices $\{1,\ldots, k,\ldots, l\}$, and $b,b',b'',\ldots$ run over the complementary set $\{l+1,\ldots,n\}$. (Note that if $l=k$, then the indices $a, a'\ldots$ run over the same set as the capital Greek indices $\Lambda,\Gamma,\ldots$, whereas $b,b',\ldots$ run over the same set as lower-case Greek indices $\alpha,\beta,\ldots$.) Let $z_0:=((x^i)_0,(p^A)_0,(p^A_i)_0)$ with
 $(p^A_\alpha)_0=F^A_\alpha((x^i)_0,(p^B)_0,(p^B_\Lambda)_0)$ and assume that $E_0\subset T_{z_0}\bR^{n+m+nm}$ is an $l$-dimensional integral element whose basis is of the form
\begin{equation}\Label{E0basis}
\frac{\partial}{\partial x^a}+(c_a^A)_0\frac{\partial}{\partial p^A} + (c_{a j}^A)_0\frac{\partial}{\partial p_j^A},\quad a=1,\ldots, l.
\end{equation}
For $z:=(x^i,p^A,p^A_i)$, let $G_l(T_z\bR^{n+m+nm})$, for $k\leq l\leq n+m+nm$, denote the Grassmannian of all $l$-dimensional subspaces of $T_z\bR^{n+m+nm}$.
Note that $(z_0,E_0)$ is a point in the Grassmannian bundle
\begin{equation}
X_l:=\big\{(z,E)\colon z\in \bR^{n+m+nm},\ E\in G_l(T_z\bR^{n+m+nm})\big\}.
\end{equation}
For any $(z,E)\in X_l$, sufficiently close to $(z_0,E_0)$, there is a (unique) basis for $E$ of the form
\begin{equation}\Label{ea}
\tilde e_a:=\frac{\partial}{\partial x^a}+c^b_a\frac{\partial}{\partial x^b}+c_a^A\frac{\partial}{\partial p^A} + c_{a j}^A\frac{\partial}{\partial p_j^A},\quad a=1,\ldots, l,
\end{equation}
where $\tilde e_a\in T_{z}\bR^{n+m+nm}$ and $(z,c^b_a,c_a^A,c_{a j}^A)$ is close to $(z_0,0,(c_a^A)_0,(c_{a j}^A)_0)$.
The map $(z,E)\mapsto (z, c^b_a,c_a^A,c_{a j}^A)$ forms a local coordinate system for $X_l$ near $(z_0,E_0)$. Consider the family $\mathcal F$ of germs at $(z_0,0,(c_a^A)_0,(c_{a j}^A)_0)$ of real-analytic functions given by
\begin{equation}\Label{familyF}
\begin{aligned}
f^A_\alpha& :=p^A_\alpha-F^A_\alpha,\\
g^A_a &:=\eta^A(\tilde e_a) =c^A_a-p^A_a-p^A_bc^b_a,\\
g^A_{\alpha a}&:=\omega^A_\alpha(\tilde e_a)
= c^A_{a\alpha}-F^A_{\alpha a} - F^A_{\alpha b}c^b_a - F^A_{\alpha B} c^B_a-F^{A \Gamma}_{\alpha B}c^B_{a\Gamma},\\
h^A_{a a'}&:=\Omega^A(\tilde e_a,\tilde e_{a'}) =c^A_{aa'} +  c^A_{ab}c^b_{a'}
-(c^A_{a'a} +  c^A_{a'b}c^b_{a}),
\end{aligned}
\end{equation}
where the $F^A_\alpha$ and their derivatives are evaluated at the
point $(x^i,p^B,p^B_\Lambda)$, and $\eta^a$, $\omega^A_\alpha$, and
$\Omega^A$ are the differential forms given by \eqref{EDS}. Observe
that, by definition, an $l$-dimensional subspace $E\subset
T_z\bR^{n+m+nm}$, with $(z,E)$ sufficiently close to $(z_0,E_0)$, is
an integral element if and only if its local coordinate $(z,
c^b_a,c_a^A,c_{a j}^A)$ is a zero of the family $\mathcal F$.

Recall
(see \cite{BCGGG}) that the integral element $E_0$ of an exterior differential system is said to be
{\it K\"ahler ordinary} if $(z_0,E_0)$ is an ordinary zero of the family
of the functions of $(z,E)$ obtained by evaluating $l$-forms from the system
on the basis of $E$ chosen as above.
In fact, one can use any basis of $E$ depending on $(z,E)$
in a real-analytic fashion.
Thus,
in case of our system $\mathcal I$, the element $E_0$ is {\it K\"ahler ordinary}
if $(z_0,E_0)$ is an ordinary zero of the family
$\mathcal F$ given by \eqref{familyF}.
Further recall that an integral element $E_0$ is said to be {\it
K\"ahler regular} if it is K\"ahler ordinary and if the dimension of
the polar spaces $H(E)$ (defined by \eqref{polar0}) is constant for
all integral elements $E\subset T_z\bR^{n+m+nm}$ (of dimension $l$)
with $(z,E)$ sufficiently close to $(z_0,E_0)$.
 We have the
following theorem, in which the notation introduced above is used.

\begin{Thm}\Label{EKO} Assume that property {\rm (ii)} of Theorem $\ref{main1}$ holds and let $l\geq k$. If $E\subset T_z\bR^{n+m+nm}$ is an an $l$-dimensional integral element with $(z,E)$ sufficiently close to $(z_0,E_0)$ in  $X_l$, then $E$ is K\"ahler regular.
\end{Thm}

For the proof of Theorem \ref{EKO}, we shall need the following lemma, whose proof is elementary and left to the reader.

\begin{Lem}\Label{EKOlem}
Let $V_{ij}(s)$, $b_i(s)$, for $i=1,\ldots, M$ and $j=1,\ldots, N_2$, be real-analytic functions in a neighborhood $W$ of $s_0\in\bR^{N_1}$. Denote by $\mathcal G$ the collection of real-analytic functions   in $W\times \bR^{N_2}$, affine in $t=(t^1,\ldots, t^{N_2})$, given by $L_i(s,t):=V_{ij}(s)t^j-b_i(s)$ for $i=1,\ldots, M$. Assume that all functions in $\mathcal G$ vanish at $(s_0,t_0)$ for some $t_0\in\bR^{N_2}$ and let $A:=\{(s,t)\in W\times \bR^{N_2}\colon f(s,t)=0,\ \forall f\in\mathcal G\}$. If there exists an integer $\kappa\geq 0$ such that, for all $s_1\in W$, $\dim \{t:(s_1,t)\in A\}=\kappa$, then  $(s_0,t_0)$ is an ordinary zero of $\mathcal G$.
\end{Lem}

\begin{proof}[Proof of Theorem $\ref{EKO}$]  If we can show that $E$ is
K\"ahler ordinary, then it follows immediately from Corollary \ref{HE} and Proposition \ref{EinTS} that $E$ is K\"ahler regular. Moreover, by definition, if $E_0$ is K\"ahler ordinary, then $E$ is also K\"ahler ordinary when $(z,E)$ is sufficiently close to $(z_0,E_0)$. Thus, it suffices to show that $E_0$ is K\"ahler ordinary. For this, we shall make use of Lemma \ref{EKOlem} with $\mathcal G$ being the family $\mathcal F$ given in the local coordinates $(z, c^b_a,c_a^A,c_{a j}^A)$ by \eqref{familyF}. Recall that $z=(x^i,p^A,p^A_i)$.
We split the local coordinates $(z, c^b_a,c_a^A,c_{a j}^A)$ into $(s,t)$, with $s=(x^i,p^A,p^A_\Lambda,c^b_a)$ and $t=(p^A_\alpha, c^A_a, c_{a j}^A)$, and take $(s_0,t_0)$ to be the local coordinates of the point $(z_0,E_0)\in X_l$. Observe that the functions in $\mathcal F$ are of the form required by Lemma \ref{EKOlem}, i.e. affine in $t$. Thus, to prove Theorem \ref{EKO} it suffices to show that the dimension of the affine planes $\{t:(s_1,t)\in A\} $ is constant in $s_1$, for $s_1$ close to $s_0$. Let us fix $s_1=((x^i)_1,(p^A)_1,(p^A_\Lambda)_1,c^b_a)$ close to $s_0$. As in the proof of Proposition \ref{EinTS}, consider the Cauchy problem \eqref{CP}
with \eqref{initial} satisfied, where $((x^i)_0,(p^A)_0,(p^A_\Lambda)_0)$ is replaced by $((x^i)_1,(p^A)_1,(p^A_\Lambda)_1)$. By assumption, this Cauchy problem has an approximate solution $u^A(x^i)$ at $x_1:=(x^i)_1$.
  By Proposition \ref{eqapprox}, the derivatives of $u^A$ at $x_1$ satisfy \eqref{altappsol}. The jet graph of $u^A$ is parametrized by the map $\Psi$ given by \eqref{Psi}. A straightforward computation, as in the proof of Proposition \ref{EinTS}, shows that
\begin{equation}\Label{e}
\begin{aligned}
e_\Lambda:= &\Psi_*\left (\partial/\partial x^\Lambda\right )=\frac{\partial}{\partial x^\Lambda}+
u_\Lambda^A\frac{\partial}{\partial p^A} + u_{\Lambda \Gamma}^A\frac{\partial}{\partial p_\Gamma^A}+(F^A_{\alpha\Lambda}+F^A_{\alpha B}u^B_\Lambda+F^{A\Gamma}_{\alpha B}u^B_{\Lambda\Gamma})\frac{\partial}{\partial p_\alpha^A},\\
e_\alpha:= &\Psi_*\left(\partial/\partial x^\alpha\right)=\frac{\partial }{\partial x^\alpha}+u^A_\alpha
\frac{\partial }{\partial p^A}+u^A_{\alpha\Lambda}\frac{\partial }{\partial p^A_\Lambda}+u^A_{\alpha\beta}\frac{\partial }{\partial p^A_\beta},
\end{aligned}
\end{equation}
where the $\partial/\partial x^i$ are tangent vectors at the point $x_1\in \bR^n$.
Recall that the coefficients $u^A_i$, $u^A_{ij}$ appearing in \eqref{e} are evaluated at $x_1$ and hence, in particular, $u^A_\Lambda=(p^A_\Lambda)_1$. Moreover, these coefficients satisfy \eqref{altappsol}, and hence we can think of them (and also then of the $e_i$) as functions of  $a^A_{\Lambda\Gamma}$. Let $z_1$ denote the point
$$((x^i)_1,(p^A)_1,(p^A_\Lambda)_1, F^A_\alpha((x^i)_1,(p^B)_1,(p^B_\Gamma)_1))\in \bR^{n+m+nm}.$$
We observe that each $e_i\in T_{z_1}\bR^{n+m+nm}$ in \eqref{e} is an affine function of $(a^A_{\Lambda\Gamma})\in \bR^{mk(k+1)/2}$, i.e.\  $e_i=e_i(a^A_{\Lambda\Gamma})$, and that the functions $e_\Omega=e_\Omega(a^A_{\Lambda\Gamma})$ are injective. Consider the map sending $(a^A_{\Lambda\Gamma})\in \bR^{mk(k+1)/2}$ to the point $(z_1,E)\in X_l$, where $E\subset T_{z_1}\bR^{n+m+nm}$ is spanned by the $l$ vectors $\tilde e_a:=e_a+c^b_ae_b$. Observe that $E$ is a subspace of $T_{z_1} S_{u^A}$, where $S_{u^A}$ is the jet graph of the approximate solution $u^A$, and hence is an integral element of $\mathcal I$. Using the local coordinates $(s,t)$ on $X_l$ introduced above, the mapping $(a^A_{\Lambda\Gamma})\mapsto (z_1,E)$ satisfies $s=s_1$ and its $t$-component, $t=t(a^A_{\Lambda\Gamma})$, is affine. Moreover, as the reader can verify, the map $t=t(a^A_{\Lambda\Gamma})$ is injective when the $c^b_a$ are sufficiently small. By Proposition \ref{EinTS}, any $l$-dimensional integral element $E\subset T_{z_1}\bR^{n+m+nm}$ with $(z_1,E)$ close to $(z_0,E_0)$, is contained in $T_{z_1} S_{u^A}$, where $S_{u^A}$ is the jet graph through $z_1$ of an
approximate solution $u^A$ at $x_1$. Consequently, for each $s_1$ close to $s_0$, the space of all $l$-dimensional integral elements $(z,E)$ with coordinates $(s_1,t)$ is parametrized by the affine map $(a^A_{\Lambda,\Gamma})\mapsto (s_1,t(a^A_{\Lambda\Gamma}))$ and hence its dimension is $\kappa=mk(k+1)/2$. By Lemma \ref{EKOlem}, $(s_0,t_0)$ is an ordinary zero for $\mathcal F$, which completes the proof of Proposition \ref{EKO}.
\end{proof}

\section{Proof of {\rm (ii)}$\implies$  {\rm (i)} in Theorem
\ref{main1}}\Label{iiitoi}

We note that it suffices to prove that the Cauchy problem \eqref{CP} is solvable at
$((x^i)_0,(p^B)_0,(p^B_\Lambda)_0)$.  We let $S_0\subset \bR^{n+m+nm}$ be the submanifold
through $z_0:=((x^i)_0,(p^B)_0,(p^B_\Lambda)_0,(p^A_\alpha)_0)$, with $(p^A_\alpha)_0:=
F^A_\alpha((x^i)_0,(p^B)_0,(p^B_\Lambda)_0)$,  parametrized by $\Phi\colon V\to \bR^{n+m+nm}$,
where $V$ is a sufficiently small open neighborhood of $(x^\Lambda)_0$ in $\bR^k$ and
\begin{equation}\Label{Phi2}
\Phi(x^\Lambda):=\big((x^\Lambda,x^\alpha_0),a^A(x^\Lambda),(a^A_\Gamma(x^\Lambda),
F^A_\alpha((x^\Lambda,x^\alpha_0),a^B(x^\Lambda),a^B_\Gamma(x^\Lambda)))\big),
\end{equation}
where $a^A(x^\Lambda)$ is the Cauchy data in \eqref{CP} satisfying
\eqref{initial}. Since the forms \eqref{E} vanish when restricted to $S_0$, the latter  is an
integral submanifold for the differential ideal $\mathcal I$ introduced
in Section \ref{ext}. Let $E_0$ be the tangent space $T_{z_0}S_0$,
and observe that $E_0$ is an  integral element of dimension $l=k$ with
a basis of the form \eqref{E0basis}. By Theorem \ref{EKO}, the integral element
$E:=T_zS_0$ is K\"ahler regular for every $z\in S_0$ with $z$
sufficiently close to $z_0$, and $\dim H(E)=n$ (by Corollary
\ref{HE} and Proposition \ref{EinTS}).
Recall that an integral submanifold is called {\em K\"ahler
regular} if its tangent space at every point is K\"ahler
regular. In view of the above observation, $S_0$ is K\"ahler
regular in some neighborhood of $z_0$. Let
$$R_1:=\{(x^i,p^A,p^A_i)\in\bR^{n+m+nm}\colon
x^{k+2}=x_0^{k+2},\ldots, x^n=x_0^n\}.$$
Note that $R_1$ has
codimension $r=n-k-1$ and is transverse to $H(E)$, for all
$E=T_zS_0$ for $z\in S_0$ near $z_0$, by Corollary \ref{HE} and
Proposition \ref{EinTS}. (Indeed, $H(E)$ is a graph over $T_x\bR^n$, being  a tangent space of the jet graph of an approximate solution, and the variables $(p^A,p^A_i)$ are free in $R_1$.) By the Cartan-K\"ahler Theorem (Theorem 2.2
in \cite{BCGGG}), there exists a real-analytic integral submanifold
$S_1$ of dimension $l=k+1$ through $z_0$ such that $S_0\subset
S_1\subset R_1$. By Proposition \ref{EinTS}, the tangent space
$T_{z_0} S_1$ is contained in $T_{z_0} S_{u^A}$ for the jet graph
$S_{u^A}$ of some approximate solution $u^A(x^i)$ at $x_0$. Since the
projection $\pi\colon \bR^{n+m+nm}\to \bR^n$ restricted to $S_{u^A}$
is a diffeomorphism, we conclude that the projection $\pi$
restricted to $S_1$ near $z_0$ is a diffeomorphism onto a
$(k+1)$-dimensional submanifold through $x_0$ in $\bR^n$. Since
$S_1$ is also contained in $R_1$, we conclude that $T_{z_0}S_1$ has
a basis of the form \eqref{E0basis}. Hence, it follows that $S_1$ is
K\"ahler regular near $z_0$ by Theorem \ref{EKO}, and we have $\dim
H(T_zS_1)=n$ for every $z\in S_1$ near $z_0$ (again by Corollary \ref{HE} and Proposition
\ref{EinTS}). Repeating the argument with $S_0$ replaced by $S_1$
and $R_1$ replaced by
$$R_2:=\{(x^i,p^A,p^A_i)\in\bR^{n+m+nm}\colon
x^{k+3}=x_0^{k+3},\ldots, x^n=x_0^n\},$$ we obtain a real-analytic
integral submanifold $S_2$ of dimension $l=k+2$ through $z_0$  such
that $S_0\subset S_1\subset S_2\subset R_2$. Again, $T_{z_0}S_2$ has
a basis of the form \eqref{E0basis}, $S_2$ is K\"ahler regular near
$z_0$, and we have $\dim H(T_zS_2)=n$ for $z\in S_2$ near $z_0$.  We continue this process,
inductively producing integral submanifolds whose dimensions increase
by one at each step. The process will end after $n-k$ steps with an
$n$-dimensional real-analytic integral submanifold $S:=S_{n-k}$ through
$z_0$ with the property that the projection $\pi\colon
\bR^{n+m+nm}\to \bR^n$ restricted to $S$ is a diffeomorphism, i.e.
$dx_1\wedge\ldots\wedge dx_n\neq 0$ on $S$. As observed in the
beginning of Section \ref{ext}, this produces a real-analytic
solution to the Cauchy problem \eqref{CP}. The uniqueness of the
solution follows easily from a standard power series argument. The
details are left to the reader. This completes the proof of the
implication (ii)$\implies$(i) in Theorem \ref{main1}.\qed

\section{A system of Monge-Amp\`ere type}\Label{Monge}

Let $u(x^i)$ be a real-valued function defined near the origin in
$\bR^n$. Denote by $Hu(x^i)$ the Hessian of $u$, i.e.\ the symmetric
$n\times n$-matrix of second order derivatives $(u_{ij})_{i,j=1}^n$.
As an application of Theorem \ref{main1}, we shall give a
construction of all real-analytic, real-valued functions $u$ defined
near $0$ in $\bR^n$ such that the rank of $Hu$ is identically one.
We shall restrict our attention to those $u$ for which
 \begin{equation}\Label{Delta}
 \delta (u):=u_{11}
 \end{equation}
 satisfies $\delta(u)(0)\neq 0$.
 Clearly every $u$ with nonvanishing $Hu$ at $0$ can be put into this form
 by a linear change of coordinates in $\bR^n$.
 Recall, from the introduction, the differential operators $\Delta_{\alpha\beta}$ given by the
 $2\times 2$-minors
\begin{equation}\Label{Deltaab}
 \Delta_{\alpha\beta}(u):=\det
 \begin{pmatrix}
 u_{11}  & u_{1\beta}\\
 u_{\alpha 1} & u_{\alpha\beta}
 \end{pmatrix},
 \end{equation}
for $\alpha$, $\beta\in\{2,\ldots, n\}$. (We will continue to use
the notation and conventions established in previous sections with
$k=1$.) Observe that
$\Delta_{\alpha\beta}(u)=\Delta_{\beta\alpha}(u)$.
The following statement is a consequence of a standard linear algebra argument:

\begin{Pro}\Label{rankk} Let $u$ be a real-valued, real-analytic function defined near $0$ in $\bR^n$. Assume that $\delta(u)(0)\neq0$. Then, the rank of $Hu$ is identically one near $0$ if and only if $\Delta_{\alpha\beta}(u)=0$ for all $\alpha$, $\beta\in\{2,\ldots, n\}$.
\end{Pro}

To find all real-analytic $u$ such that the rank of $Hu$ is
identically one and $\delta(u)(0)\neq 0$, we shall consider the
overdetermined Cauchy problem
\begin{equation}\Label{MA0}
\left\{
\begin{aligned}
 & \Delta_{\alpha\beta}(u)=0,\\ &u(x^1,0)=a(x^1),\quad u_{\alpha}(x^1,0)=a_\alpha(x^1),
\end{aligned}
\right . \quad \alpha,\beta=2,\ldots, n,
\end{equation}
where the $n$ functions $a(x^1)$, $a_\alpha(x^1)$ are real-valued and
real-analytic near $0$ in $\bR$. In order to ensure that
$\delta(u)(0)\neq 0$ we shall require $a''(0)\neq 0$. We shall in
fact consider the more general Cauchy problem, as in Theorem \ref{MAthm},
 \begin{equation}\Label{MA}
\left\{
\begin{aligned}
 & \Delta_{\alpha\beta}(u)=f_{\alpha\beta}(x^j,u_{11}),\\ &u(x^1,x^\alpha_0)=a(x^1),\quad u_{\alpha}(x^1,x^\alpha_0)=a_\alpha(x^1),
\end{aligned}
\right . \quad \alpha,\beta=2,\ldots, n,
\end{equation}
where the $f_{\alpha\beta}(x^j,t)$, $(x^j)_0$, $a(x^1)$, and $a_\alpha(x^1)$ are as in Theorem \ref{MAthm}. Naturally, we have to impose the
symmetry condition $f_{\alpha\beta}=f_{\beta\alpha}$. We now give the proof
of Theorem \ref{MAthm}.

\begin{proof}[Proof of Theorem $\ref{MAthm}$] We shall reduce the Cauchy problem \eqref{MA} to a Cauchy problem for a
system of first order partial differential equations of the form \eqref{CP},
and verify that
the compatibility conditions of (iii) in Theorem
\ref{main1}  are equivalent to the $f_{\alpha\beta}$  being of the form $f_{\a\b}(x^j,t)=g_{\a\b}(x^j)t$ with $g_{\a\b}$ satisfying \eqref{symder0}.
The conclusion of Theorem \ref{MAthm} then
follows from Theorem \ref{main1}. In the notation used in Theorem
\ref{main1}, we let $m=n+1$ and, for a real-valued $u$, we define
$u^A$, $1\leq A\leq n+1$, as follows
\begin{equation}\Label{uA}
u^i:=u_i=\frac{\partial u}{\partial x^i},\ 1\leq i\leq n,\quad u^{n+1}:=u.
\end{equation}
We define, for a vector valued function $u^i$ with $1\leq i\leq n$, the first order differential operators
\begin{equation}\Label{Deltaab2}
 \Delta_{\alpha}^{\beta}(u^i):=\det
 \begin{pmatrix}
 u_{1}^{1} & u_{1}^{\beta}\\
 u_{\alpha }^{1} &  u_{\alpha}^{\beta}
 \end{pmatrix},
 \end{equation}
for $\alpha$, $\beta\in\{2,\ldots, n\}$. By construction, a solution $u$ to the second order Cauchy problem \eqref{MA} yields a solution $u^A$, via \eqref{uA}, to
the first order Cauchy problem:
\begin{equation}\Label{MA3}
\left\{
\begin{aligned}
 & u^1_\alpha=u^\alpha_1,\\
 & \Delta^\beta_\alpha(u^i)=f_{\alpha\beta}(x^j,u^1_1),\\
 &u^{n+1}_\alpha=u^\alpha,\\
 &u^1(x^1,x^\a_0)=a'(x^1),\ u^{\alpha}(x^1,x^\a_0)=a_\alpha(x^1),\\ & u^{n+1}(x^1,x^\a_0)=a(x^1),
\end{aligned}
\right.
\quad\alpha,\beta=2,\ldots, n.
\end{equation}
Note that we did not include the equation $u^{n+1}_1=u^1$
which does not contain a derivative transversal to the line $\{(x^1,x^\a_0)\}$.
However, it actually follows from other equations.
In fact, we have the converse:

\begin{Lem} \Label{inverse} If $u^A$, $1\leq A\leq n+1$, is a solution to the Cauchy problem \eqref{MA3}, then $u:=u^{n+1}$ is a
solution to the Cauchy problem \eqref{MA} (as a germ at $(x^j)_0$).
\end{Lem}

\begin{proof}  Clearly, $u:=u^{n+1}$ satisfies the boundary condition  (i.e.\ the last line)
in \eqref{MA}. Thus, to prove the lemma, we must show that
$\Delta_{\alpha\beta}(u)= f_{\alpha\beta}(x^j, u_{11})$. To do this it suffices to
show that
\begin{equation}\Label{deltaids}
\begin{aligned}
& u_{\alpha\beta}=u^\alpha_\beta,\\
& u_{1 \beta}=u^\beta_1,\ u_{\alpha 1}=u^1_\alpha,\\
& u_{11}=u_1^1,
\end{aligned}
\end{equation}
since this implies
$\Delta_{\alpha\beta}(u)=\Delta^\beta_\alpha(u^i)= f_{\alpha\beta}(x^j, u_{11})$.
The first two lines in \eqref{deltaids} follow directly from the
differential equations in \eqref{MA3} (using the symmetry of $u_{\a\b}$ in $\a$ and $\b$).
To check the remaining
identity in \eqref{deltaids}, we observe that the initial data
conditions in \eqref{MA3} imply that it holds when $x^\alpha=x^\a_0$. On
the other hand, by differentiating the first and third lines in
\eqref{MA3}, we obtain
$$
u_{\alpha 11}=u^{n+1}_{\alpha11}=u^\alpha_{11}=u^1_{\alpha 1},
$$
and, hence, $(u_{11}-u^1_1)_\alpha=0$. Since
$u_{11}=u^1_1$ when $x^\alpha=x^\a_0$, we conclude, by a uniqueness argument, that the identities in the last line of \eqref{deltaids} hold. This completes the proof of the lemma.
\end{proof}

We now observe that the Cauchy problem \eqref{MA3} decouples into the two Cauchy problems
\begin{equation}\Label{MA3A}
\left\{
\begin{aligned}
 & u^1_\alpha=u^\alpha_1,\\
 & \Delta^\beta_\alpha(u^i)=f^\beta_\alpha(x^j,u^1_1),\\
 &u^i(x^1,x^\a_0)=a^i(x^1),
\end{aligned}
\right.
\quad\alpha,\beta=2,\ldots, n,
\end{equation}
 and
 \begin{equation}\Label{MA3B}
\left\{
\begin{aligned}
 &u^{n+1}_\alpha=u^\alpha,\\
  & u^{n+1}(x^1,x^\a_0)=a(x^1),
\end{aligned}
\right.
\quad\alpha,\beta=2,\ldots, n,
\end{equation}
where the Cauchy data $a^i$ in \eqref{MA3A} are taken to be
$a^1:=a'$ and $a^\alpha:= a_\alpha$, and the right hand side
$f^\beta_\alpha:=f_{\alpha\beta}$.

\begin{Pro}\Label{MA3Asol} Let
$f_{\alpha}^{\beta}(x^i,t)$, for $2\leq \alpha,\beta\leq n$, with $n>2$,
be real-analytic functions in a connected open subset
$U\times V\subset\bR^{n}\times(\bR\setminus\{0\})$.
Then the following two conditions are equivalent:

\medskip
\noindent {\rm (i)} For any choice of real-analytic Cauchy data $a^i(x^1)$ near $x_0$ in $U$ such that $(a^1)'(x^1_0)\in V$, there exists a unique real-analytic solution $u^j(x^i)$ near $x_0$ to the
Cauchy problem \eqref{MA3A}.

\medskip
\noindent
{\rm (ii)} The functions $f_{\alpha}^{\beta}$ are of the form
$f_{\a}^\b(x^i,t)=g_{\a}^\b(x^i)t$ with $g_{\a}^\b$ satisfying
\begin{equation}\Label{symder1}
 g_{\alpha 1}^{\beta}=0,\quad
g_{\alpha\beta}^{\gamma}=g_{\beta \alpha}^{\gamma},\quad
\alpha,\beta,\gamma=2,\ldots, n.
\end{equation}
\medskip

Moreover, if $f^\beta_\alpha=f^\alpha_\beta$, then any solution $u^j(x^i)$ to \eqref{MA3A} satisfies $u_\beta^\alpha=u^\beta_\alpha$ for $2\leq \alpha,\beta\leq n$.
\end{Pro}

\begin{proof} Let us rewrite the differential equations in
\eqref{MA3A} explicitly in the form of \eqref{CP}:
\begin{equation}\Label{MA4A}
\left\{
\begin{aligned}
 & u^1_\alpha=u^\alpha_1,\\
 & u^\beta_\alpha=\frac{u^\beta_1u^\alpha_1+
f^\beta_\alpha}{u^1_1},
\end{aligned}
\right. \quad\alpha,\beta=2,\ldots, n.
\end{equation}
In the second line, we have used the first line to replace
$u^1_\alpha$ by $u^\alpha_1$. We remark that, in the notation of
Theorem \ref{main1}, we have here $m=n$ and $k=1$, so that the indices
$A$, $B$, etc.\ run over the set $\{1,\ldots, n\}$ (i.e.\ the same
index set as $i$, $j$, etc.) and $\Lambda$, $\Gamma$, etc.\, are all
$1$, and the right-hand side functions are given by
\begin{equation}
F^1_\alpha (x^i, p^B, p^B_1) = p_1^\alpha, \quad
F^\beta_\alpha(x^i, p^B, p^B_1)=\frac{p^\beta_1p^\alpha_1+
f^\beta_\alpha(x^i,p^1_1)}{p^1_1}.
\end{equation}
We compute the derivatives of $F^A_\a$ that appear in \eqref{comp10} and
\eqref{comp20} and that are not $0$:
\begin{equation}\Label{Ders}
\left\{
\begin{aligned}
F^\gamma_{\alpha i} & =\frac{f^{\gamma}_{\alpha i}(x^j,p^1_1)}{p^1_1},\\
F^{11}_{\alpha \alpha} & =1,\\
F^{\beta 1}_{\alpha 1} & = -\frac{p^\alpha_1 p^\beta_1 + f^\b_\a}{(p^1_1)^2}
+\frac{f^{\b1}_{\a1}}{p^1_1},\\
F^{\beta 1}_{\alpha \gamma} & =
\frac{\delta_{\alpha\gamma}p^{\beta}_1+\delta_{\beta\gamma}p^\alpha_1}{p^1_1},
\end{aligned}
\right.
\end{equation}
where $\delta_{\alpha\gamma}=\delta_\alpha^\gamma$ denotes the
Kronecker delta symbol. We remark
that the compatibility conditions \eqref{comp10} and \eqref{comp20}
simply say that the functions in their left hand sides (or, equivalently, the right
hand sides) are symmetric in $\alpha$ and $\beta$ for any $A\in
\{1,\ldots, n\}$ in the case of \eqref{comp10}, and for any $A, C\in
\{1,\ldots, n\}$ and $\Gamma=\L=1$ in the case of \eqref{comp20}. We compute the left
hand side of \eqref{comp10} with $A=1$ to obtain
\begin{equation}\Label{1}
 \Phi^1_{\a\b}=0+0+F^{11}_{\alpha B}(F^B_{\beta 1} + 0)
 =F^\alpha_{\beta 1} =\frac{f^\alpha_{\beta 1}}{p^1_1}.
\end{equation}
 Next, we compute
the left hand side of \eqref{comp10} with $A=\gamma$:
\begin{equation}\Label{2}
\Phi^\gamma_{\a\b}=F^\gamma_{\alpha\beta}+0+F^{\gamma 1}_{\alpha B}(F^B_{\beta
1}+0)=\frac{f^\gamma_{\a\b}}{p^1_1}+
\frac{
\delta_{\alpha\epsilon}p^{\gamma}_1+\delta_{\gamma\epsilon}p^\alpha_1}
{(p^1_1)^2}f^\epsilon_{\b1}=
\frac{f^\gamma_{\a\b}}{p^1_1}+\frac{p^\gamma_1f^\alpha_{\beta 1} + p^\alpha_1f^\gamma_{\beta 1}}{(p^1_1)^2}.
\end{equation}
  We now compute the
left hand side of \eqref{comp20} with $A=C=1$:
\begin{equation}\Label{3}
\Psi^{111}_{\a\b1}
= 2F^{11}_{\alpha B} F^{B 1}_{\beta 1}
=-2\frac{p^\alpha_1 p^\beta_1 + f^\b_\a}
{(p^1_1)^2} + 2\frac{f^{\b1}_{\a1}}{p^1_1}.
\end{equation}
The left hand side of \eqref{comp20} with $A=1$, $C=\gamma$ is given
by
\begin{equation}\Label{4}
\Psi^{111}_{\a\b \gamma}
=2F^{11}_{\alpha B} F^{B 1}_{\beta
\gamma}
=2\frac{\delta_{\alpha\gamma}p^{\beta}_1+\delta_{\beta\gamma}p^\alpha_1
}{p^1_1},
\end{equation}
and
with $A=\gamma$, $C=1$ by
\begin{multline}\Label{5}
\begin{aligned}
\Psi^{\gamma11}_{\a\b 1}
=2F^{\gamma 1}_{\alpha B} F^{B 1}_{\beta 1}
=&2(\delta_{\alpha\epsilon}p^\gamma_1+\delta_{\gamma\epsilon}p^\alpha_1)
\Big(-\frac{p^\epsilon_1 p^\b_1 +  f_\b^\epsilon}{(p^1_1)^3}
+
\frac{f^{\epsilon1}_{\b1}}{(p^1_1)^2}\Big)\\
=&-2\frac{2p^\alpha_1 p^\beta_1 p^\gamma_1
+ p^\gamma_1 f^\a_\b + p^\a_1 f^\gamma_\b}{(p^1_1)^3}
+ 2\frac{p^\gamma_1 f^{\a1}_{\b1} + p^\a_1 f^{\gamma1}_{\b1}}{(p^1_1)^2},
\end{aligned}
\end{multline}
and finally with $A=\gamma$ and $C=\lambda$ by
\begin{equation}\Label{6}
\begin{aligned}
\Psi^{\gamma11}_{\a\b \lambda}=
2F^{\gamma 1}_{\alpha B} F^{B 1}_{\beta \lambda} &
= 2F^{\gamma
1}_{\alpha 1} F^{1 1}_{\beta \lambda} + 2F^{\gamma 1}_{\alpha
\epsilon} F^{\epsilon 1}_{\beta \lambda} \\&
=2\Big(-\frac{p^\gamma_1 p^\alpha_1 + f^\gamma_\a}{(p^1_1)^2}
+\frac{f^{\gamma1}_{\a1}}{p^1_1}\Big)
\delta_{\beta\lambda}
+2\frac{(\delta_{\alpha\epsilon}p^\gamma_1+
\delta_{\gamma\epsilon}p^\alpha_1)
 (\delta_{\beta\lambda}p^\epsilon_1+
\delta_{\lambda\epsilon}p^\beta_1)} {(p^1_1)^2} \\ &=
2\frac{\delta_{\alpha\lambda} p^\gamma_1 p^\beta_1+
\delta_{\beta\lambda} p^\alpha_1 p^\gamma_1 + \delta_{\gamma\lambda}
p^\beta_1 p^\alpha_1}{(p^1_1)^2}
-2\Big(\frac{f^\gamma_\a}{(p^1_1)^2}-\frac{f^{\gamma1}_{\a1}}{p^1_1}\Big)\delta_{\b\lambda}.
\end{aligned}
\end{equation}
Since $n>2$,
the expressions in \eqref{6} are symmetric in $\a$ and $\b$
if and only if
$$\frac{f^\gamma_\a}{(p^1_1)^2}-\frac{f^{\gamma1}_{\a1}}{p^1_1}=0,
\quad \a,\gamma=2,\ldots,n,$$
or, equivalently,
$$
f^{\gamma1}_{\a1}=\frac{f^\gamma_\a}{p^1_1}.
$$
The latter is an ordinary differential equation (with parameters $x^i$)
whose solutions are of the form
\begin{equation}\Label{form}
f^\gamma_\a(x^i,t)=g^\gamma_\a(x^i)t.
\end{equation}
In view of \eqref{form}, the expressions \eqref{1}-\eqref{6} become
\begin{equation}\Label{1'}
 \Phi^1_{\a\b}=g^\alpha_{\beta 1},\quad
\Phi^\gamma_{\a\b}
=g^\gamma_{\alpha\beta} +\frac{p^\gamma_1g^\alpha_{\beta 1} + p^\alpha_1g^\gamma_{\beta 1}}{p^1_1},
\end{equation}
\begin{equation}\Label{3'}
\Psi^{111}_{\a\b1}=-2\frac{p^\alpha_1
p^\beta_1}{(p^1_1)^2},
\quad
\Psi^{111}_{\a\b \gamma}
=2\frac{\delta_{\alpha\gamma}p^{\beta}_1+\delta_{\beta\gamma}p^\alpha_1}{p^1_1},
\quad
\Psi^{\gamma11}_{\a\b 1}=
-4\frac{p^\alpha_1 p^\beta_1 p^\gamma_1}{(p^1_1)^3},
\end{equation}
\begin{equation}\Label{6'}
\begin{aligned}
\Psi^{\gamma11}_{\a\b \lambda}=
2\frac{\delta_{\alpha\lambda} p^\gamma_1 p^\beta_1+
\delta_{\beta\lambda} p^\alpha_1 p^\gamma_1 + \delta_{\gamma\lambda}
p^\beta_1 p^\alpha_1}{(p^1_1)^2}.
\end{aligned}
\end{equation}
We note that all the expressions \eqref{3'}-\eqref{6'} are symmetric
in $\alpha$ and $\beta$. If we take $p^\gamma_1=p^\alpha_1=0$ in
the second expression in
\eqref{1'}, we see that this expression is symmetric in $\alpha$ and $\beta$ if and only if $g^\gamma_{\alpha\beta}$ is. Now, we take $p^\alpha_1=0$ and $p^\gamma_1=1$ in the second expression in \eqref{1'}. Since $g^\gamma_{\alpha\beta}$ must be symmetric in $\alpha$, $\beta$, we conclude that $g^\alpha_{\beta 1}$ must be symmetric in $\alpha$, $\beta$. Finally, since  $g^\gamma_{\alpha\beta}$  and $g^\alpha_{\beta 1}$ must be symmetric in $\alpha$, $\beta$, we conclude that
$p^\alpha_1  g^\gamma_{\beta 1}$ must also be symmetric in $\alpha$, $\beta$. However, since $ g^\gamma_{\beta 1}$ is a function of $x^i$ alone, this can only happen if  $ g^\gamma_{\beta 1}=0$ (where we again use $n>2$). Hence, the implication (i) $\implies$ (ii) follows from Theorem \ref{main1}. Moreover, if the conditions in (ii) are satisfied, then all the expressions \eqref{1}-\eqref{6} are symmetric in $\alpha$, $\beta$ and, hence, the  implication (ii)$\implies$(i) also follows from Theorem \ref{main1}. If $f^\alpha_\beta=f^\beta_\alpha$, then the symmetry $u^\alpha_\beta=u^\beta_\alpha$ is immediate from the
second line in \eqref{MA4A}.
\end{proof}

We now complete the proof of Theorem \ref{MAthm}. To prove (i)$\implies$(ii), we first note that the symmetry $f_{\alpha\beta}=f_{\beta\alpha}$ obviously follows from \eqref{MA00}. We consider the first order Cauchy problem \eqref{MA3A} with $f^\beta_\alpha:=f_{\alpha\beta}$.
Assume (i) in Theorem~\ref{MAthm}. We claim that (i) of Proposition \ref{MA3Asol} holds. Let $a^j(x^1)$ be real-analytic Cauchy data near $x_0$ in $\bR$ with $(a^1)'(x_0)\neq0$, and let $a(x^1)$ be a real-analytic function near $x_0$ such that $a'(x^1)=a^1(x^1)$. Note that $a''(x_0)=(a^1)'(x_0)\neq0$. We now let $u(x^i)$ be the solution to the Cauchy problem  \eqref{MA00} with Cauchy data $a(x^1)$ and $a_\alpha(x^1):=a^\alpha(x^1)$ (that exists in view of our assumption). The gradient $u^j(x^i):=u_j(x^i)$ solves the Cauchy problem \eqref{MA3A}. Hence, by Proposition \ref{MA3Asol}, the conditions in \eqref{symder0} must also hold.

To prove the converse (ii)$\implies$(i), we let $a(x^1)$ and $a_\alpha(x^1)$ be real-analytic Cauchy data near $x_0$ in $U$ such that $a''(x^1_0)\in V$. We define $a^1(x^1):=a'(x^1)$ and $a^\alpha(x^1):=a_\alpha(x^1)$ and note that $(a^1)'(x_0)\neq 0$. We consider the Cauchy problem \eqref{MA3A} with $f^\beta_\alpha:=f_{\alpha\beta}$ and Cauchy data $a^j(x^1)$ defined above. The  conditions in \eqref{symder0} clearly imply those in \eqref{symder1}. Hence, by Proposition \ref{MA3Asol}, the Cauchy problem \eqref{MA3A} has a real-analytic solution $u^j(x^i)$ and since $f^\beta_\alpha=f^\alpha_\beta$,  we have $u^\alpha_\beta=u^\beta_\alpha$. Hence, the Cauchy problem \eqref{MA3B}, with data $a(x^1)$ as above, has a solution $u:=u^{n+1}$. By Lemma \ref{inverse}, the real-analytic function $u(x^i)$ is the solution to \eqref{MA}. This completes the proof of Theorem \ref{MAthm}.
\end{proof}

\section{Proof of Theorem \ref{Gauss}}\Label{realgeom}

In this section, we shall give the proof of Theorem \ref{Gauss}
stated in the introducion. We retain the notation established in
previous sections. We
shall use $x$ as local coordinates on $\Sigma$ near $0$. We choose as a normal to $\Sigma$ at a point $(x,u(x))$ the vector
\begin{equation}
N(x):=(u_{i}(x),-1)
\end{equation}
and, hence, the Gauss map $G\colon \Sigma\to S^{n}\subset \bR^{n+1}$ is
given by
\begin{equation}
G(x)=\frac{N(x)}{|N(x)|}.
\end{equation}
Consider the mapping $P\colon S^{n}\cap \{(x,y)\in \bR^{n+1}\colon
y<0\}\to \bR^n$ given by
\begin{equation}\Label{stereo}
P\big (x,-\sqrt{1-|x|^2}\big ):=\frac{x}{\sqrt{1-|x|^2}}.
\end{equation}
Note that the mapping $\big (x,-\sqrt{1-|x|^2}\big )\mapsto
(P(x),-1)$ is the stereographic projection from the origin to the
hyperplane $\{(x,y)\in \bR^{n+1}\colon y=-1\}$. Consequently, we
conclude that $P$ is a diffeomorphism and we have the identity
$P\circ G=(u_{i})$. Let $\Gamma\colon(-\epsilon,\epsilon)\to \bR^n$,
for $\epsilon>0$ sufficiently small, be the real-analytic
curve $\Gamma:=P\circ(\gamma|_{(-\epsilon,\epsilon)})$ through the origin in $\bR^n$. By assumption, $\Gamma$ is of the form
\begin{equation}
\Gamma(t)=(\varphi(t), a_2(t),\ldots, a_n(t)),
\end{equation}
with $\varphi(0)=0$, $\varphi'(0)=1$, $a'_\a(0)=0$.
We define $a(x^1)$ by
$$a(x^1):=\int_{0}^{x^1}\varphi(t) \, dt.$$
Note that
$a''(0)=1\neq 0$. Let $u(x)$ be the unique real-analytic
solution to the Cauchy problem \eqref{MA}, with
$f_{\alpha\beta}\equiv 0$ and $x_0=0$, given by Theorem \ref{MAthm}. Now, by
construction, $u_x(x^1,0)=\Gamma(x^1)$.
Moreover, the Jacobian matrix of the gradient map
$x\mapsto u_x(x)$ equals the Hessian of $u(x)$ and hence, by
Proposition \ref{rankk}, the rank of this Jacobian matrix is one
near $0$. This completes
the proof of Theorem \ref{Gauss}.\qed

\section{Proof of Proposition \ref{Tube}}\Label{cplxgeom}

In this section, we give the proof of Proposition \ref{Tube}. Let
$\Sigma\subset \bR^{n+1}$ satisfy the assumptions in the proposition.
We shall prove that there is an open set in the tube hypersurface
$M:=\Sigma+i\bR^{n+1}$ on which $M$ is $n$-nondegenerate. The fact that
$M$ is holomorphically nondegenerate then follows from the
connectedness of $M$ and a general result relating finite
nondegeneracy and holomorphic nondegeneracy (see e.g.\ Theorem
11.5.1 in \cite{BERbook}). Moreover, it follows from the
real-analyticity of the vector valued functions in \eqref{k-nondeg},
as $p$ varies over $M$, that $M$ is $n$-nondegenerate on a dense
open set. The fact that $\ell(M)=n$ now follows from the lower
semicontinuity of the rank of a collection of vector valued
functions. Finally, the fact that the rank of the Levi form is one
on a dense open set follows easily from the fact that the rank $r$
of the Levi form at $p\in M$ equals $R-1$, where $R$ denotes the
rank of the collection of vectors in \eqref{k-nondeg} at $p$ with $j=1$ and the fact that, on a dense open set, the rank of these vectors has to strictly increase with $j$.

After a rotation and translation in $\bR^{n+1}$ if necessary, we may
assume that $\Sigma$ is given as a graph through $0$ of the form
$y=u(x)$. Since the Gauss image $G(\Sigma)$ is assumed to be a
curve, we may also assume that the rank of the Gauss map, or
equivalently the rank of the Hessian of the graphing function $u(x)$
(see the discussion in the previous section), equals one in a
neighborhood of $0$, and that $u_{11}(0)\neq 0$. Now, the tube
hypersurface $M$ is then defined in a tube neighborhood
$U+i\bR^{n+1}$, for some open neighborhood $0\in U\subset
\bR^{n+1}$, by an equation of the form $\re w=u(\re z)$, where
$(z,w)$ are coordinates in $\bC^n\times\bC$ such that $\re w=y$ and
$\re z=x$. We can also write the equation of $M$ in the form
\begin{equation}
w=-\bar w+2i u((z+\bar z)/2).
\end{equation}
If we write $q(z,\bar z):=2iu((z+\bar z)/2)$, then it is well known
(see e.g.\ Corollary 11.2.14 in \cite{BERbook}) that $M$ is
$l$-nondegenerate at at point $(z,w)\in M$ if and only if the
collection of vectors
\begin{equation}\Label{k-nondeg2}
q_{\bar z z^I} (z,\bar z),\quad I\in \bZ_+^n,\ |I|\leq j,\
\end{equation}
spans $\bC^n$ with $j=l$ and $l$ is the smallest integer with this
property; we use here the notation $q_{\bar z}$ for the vector given
by the anti-holomorphic gradient $q_{\bar z}:=(q_{\bar z^i})\in
\bC^n$, and also multi-index notation as in \eqref{k-nondeg}, so
that $z^I:=(z^1)^{I_1}\ldots (z^n)^{I_n}$. Thus, it follows from
\eqref{k-nondeg2} that $M$ is $l$-nondegenerate along the fiber
$(x,y)+i\bR^{n+1}$, with $(x,y)\in \Sigma$, if and only if the
collection of vectors
\begin{equation}\Label{k-nondeg3}
u_{x x^I} (x),\quad I\in \bZ_+^n,\ |I|\leq j,\
\end{equation}
spans $\bR^n$ with $j=l$ and $l$ is the smallest integer with this
property. We have the following lemma.

\begin{Lem}\Label{nohyper} Let $\Sigma$ be given by $y=u(x)$. The Gauss image
$G(\Sigma)$ is not contained in a hyperplane if and only if the vectors
in \eqref{k-nondeg3} evaluated at $x=0$ span $\bR^n$ for some
integer $j$.
\end{Lem}

\begin{proof} Let $P\colon S^{n}\cap \{(x,y)\in \bR^{n+1}\colon y<0\}\to \bR^n$
denote the stereographic projection introduced in \eqref{stereo} and
recall that $u_x(x)=(P\circ G)(x)$. It is easy to see that $G(\Sigma)$
is not contained in a hyperplane in $\bR^{n+1}$ if and only if
$P(G(\Sigma))$ is not contained in a hyperplane in $\bR^n$. Moreover,
$P(G(\Sigma))$ is contained in hyperplane if and only if there is a
non-zero vector $a=(a^i)\in\bR^n$ such that $a^iu_{x^i}(x)=0$ for
all $x$ in a neighborhood of $0$. Since $u$ is real-analytic, the
latter is equivalent to $a^i u_{x^i x^I}(0)=0$ for all multi-indices
$I\in \bZ_+^n$, which clearly is equivalent to the vectors in
\eqref{k-nondeg3} evaluated at $x=0$ not spanning $\bR^n$ for any
integer $j$. This completes the proof of the lemma.
\end{proof}

Since the rank of the Hessian of $u$ is one near 0 and
$u_{11}(0)\neq 0$, we conclude by Proposition \ref{rankk} that $u$
satisfies near $0$ the differential equations
$\Delta_{\alpha\beta}(u)=0$, where $\Delta_{\alpha\beta}$ are given
by \eqref{Deltaab}. Thus, for any
$\alpha\in\{2,\ldots,n\}$,
\begin{equation}\Label{gradder}
u_{x\alpha}=\frac{u_{1\alpha}}{u_{11}}u_{x1},
\end{equation}
where $u_{xi}=(u_{1i},\ldots,u_{ni})$.
An induction, left to the reader, using the identity \eqref{gradder}
proves the following lemma.

\begin{Lem}\Label{lindep} Let $u(x)$ be a function such that the Hessian $Hu$ has rank one
at every point near
$0$ and $u_{11}(0)\neq 0$. Then, for any multi-index $I\in\bZ_+^n$,
the vector $u_{x x^I}(x)$ belongs to the span of the vectors
\begin{equation}
u_{x (x^1)^m}(x),\quad m\leq |I|.
\end{equation}
\end{Lem}

Thus, we conclude, by Lemma \ref{nohyper} and \ref{lindep} that, for
some integer $j$, the collection of vectors $u_{x(x^1)^m}(0)$, with
$m\leq j$, spans $\bR^n$. It is easily seen that this is impossible
unless the collection of vectors $u_{x(x^1)^m}(x)$, with $m\leq n$,
spans $\bR^n$ for all $x$ in a dense open set. In view of Lemma
\ref{lindep} and the remarks preceding Lemma \ref{nohyper}, this
proves that $M$ is $n$-nondegenerate on an open set and, hence,
$\ell(M)=n$ as claimed. This completes the proof of Proposition
\ref{Tube}.\qed

We conclude this paper with two remarks.

\begin{Rem}\Label{prelast}
{\rm As consequence of Lemma~\ref{nohyper}
and the discussion preceding it we conclude that for a connected real-analytic tube hypersurface $M=\Sigma+i\bR^{n+1}$, the following are equivalent:
\begin{enumerate}
\item $M$ is holomorphically nondegenerate;
\item $M$ is finitely nondegenerate at each point;
\item the Gauss image $G(\Sigma)$ is not contained in a hyperplane of $\bR^{n+1}$.
\end{enumerate}
Note that, for a general connected real-analytic hypersurface of $\bC^{n+1}$,
(2) implies (1) but not vice versa.}
\end{Rem}

\begin{Rem}\Label{lastrem} {\rm It follows from the results in the
last three sections (\ref{Monge}, \ref{realgeom}, and
\ref{cplxgeom}) that, for any $n$ real-analytic functions $a(x^1)$
and $a_\alpha(x^1)$ near $0$ in $\bR$ such that $a_{11}(0)\neq 0$,
there is a unique tube hypersurface $M$ of the form $\re w=u(\re z)$
such that the rank of the Levi form is one at every point near $0$
and such that $u(x^1,0)=a(x^1)$ and $u_\alpha(x^1,0)=a_\alpha(x^1)$.
The tube hypersurface will be finitely nondegenerate at $0$, and
hence holomorphically nondegenerate, if the collection of vectors
$$
(a^{(m)}(0),a_\alpha^{(m-1)}(0)),\quad m=1,2,\ldots,
$$
spans $\bR^n$.
 }
\end{Rem}


\begin{thebibliography}{BCG$^3$91}

\bibitem[AG04]{AG04} Akivis, M. A.; Goldberg, V. V.: {\em Differential
Geometry of Varities with Degenerate Gauss Maps.} CMS Books in
Mathematics, Springer-Verlag, New York, NY, 2004.

\bibitem[BER99]{BERbook} Baouendi, M. S.; Ebenfelt, P.; Rothschild, L. P.: {\em Real
submanifolds in complex space and their mappings.} Princeton
Mathematical Series, 47. Princeton University Press, Princeton, NJ,
1999. xii+404 pp.

\bibitem[BHR96] {BHR} Baouendi, M. S.; Huang, Xiaojun; Preiss
Rothschild, L. P.: Regularity of CR mappings between algebraic
hypersurfaces. {\em Invent. Math.}, {\bf 125}, (1996), no. 1,
13--36.

\bibitem[BCG$^3$91] {BCGGG} Bryant, R. L.; Chern, S. S.; Gardner, R. B.;
Goldschmidt, H. L.; Griffiths, P. A.: {\em Exterior differential
systems.} Mathematical Sciences Research Institute Publications, 18.
Springer-Verlag, New York, 1991. viii+475 pp.

\bibitem[C31]{C31} Cartan, E.: Sur la thŽorie des systmes en involution et ses applications ˆ la relativitŽ. (French) {\em  Bull. Soc. Math. France}, {\bf  59},  (1931), 88--118.

\bibitem[E01]{E01} Ebenfelt, P.: Uniformly Levi degenerate CR manifolds: the
5-dimensional case.  {\em Duke Math. J.}, {\bf  110},   (2001),  no.
1, 37--80. Erratum: {\em Duke Math. J.}, {\bf  131},  (2006),  no.
3, 589--591.

\bibitem[F07]{F07} Fels, G.: Locally homogeneous finitely nondegenerate
CR-manifolds.  {\em Math. Res. Lett.},  {\bf 14},  (2007),  no. 6,
893--922.

\bibitem[FK06]{FK06} Fels, G.; Kaup, W.: Classification of Levi degenerate
homogeneous CR-manifolds in dimension 5. {\em Acta Math} (to
appear). {\tt http://arxiv.org/ps/math.CV/0610375}

\bibitem[FK07]{FK07} Fels, G.; Kaup, W.: CR-manifolds of
dimension 5: a Lie algebra approach. {\em J. Reine Angew. Math.},
{\bf  604} (2007), 47--71.

\bibitem[F77]{F77} Freeman, M.: Local biholomorphic straightening of real
submanifolds.  {\em Ann. Math.}, (2), {\bf 106},   (1977), no. 2,
319--352.

\bibitem [G67]{G67} Goldschmidt, H.:
Integrability criteria for systems of nonlinear partial differential equations.
{\em J. Differential Geometry}, {\bf  1}, (1967), 269--307.


\bibitem[HN59]{HN} Hartman, P.; Nirenberg, L.:
On spherical image maps whose Jacobians do not change sign. {\em
Amer. J. Math.}, {\bf 81}, (1959), 901--920.

\bibitem[R10]{R10} Riquier, C.: {\em Les syst\`emes d'\'equations aux d\'eriv\'ees partielles.}, Gauthiers-Villars, Paris, 1910.

\bibitem[S60]{S60} Sacksteder, R.: On hypersurfaces with no negative sectional
curvatures. {\em Amer. J. Math.}, {\bf 82}, (1960), 609--630.

\bibitem[T34]{T34} Thomas, J. M.: Riquier's existence theorems. {\em  Ann. of Math.} (2), {\bf  35},  (1934),  no. 2, 306--311.

\bibitem[W95] {W95} Wu, H.: Complete developable submanifolds in real and complex
Euclidean spaces. {\em Internat. J. Math.}, {\bf 6}, (1995), no. 3,
461--489.


\end{thebibliography}
\end{document}